\documentclass[10pt,journal]{IEEEtran}
\usepackage{ifpdf}
 \ifpdf
 \else
 \fi
\usepackage{graphicx}
\graphicspath{{./figures/}}
\usepackage{amsmath,amssymb,amsfonts}
\usepackage{textcomp}
\allowdisplaybreaks[4]
\usepackage{makecell}
\usepackage{xcolor}
\usepackage{color}
\usepackage{caption}
\usepackage{booktabs}% 表格线粗细

\usepackage{algorithm}  
\usepackage{algorithmic} 
\makeatletter
\let\NAT@parse\undefined
\makeatother
\usepackage{hyperref}  %hyperref still needs to be put at the end!
\newtheorem{definition}{Definition}
\newtheorem{assumption}{Assumption}
\newtheorem{lemma}{Lemma}
\newtheorem{theorem}{Theorem}

\newtheorem{remark}{Remark}

\ifCLASSOPTIONcompsoc
  \usepackage[nocompress]{cite}
\else
  \usepackage{cite}
\fi
\ifCLASSINFOpdf
\else
\fi
\ifCLASSOPTIONcompsoc
  \usepackage[caption=false,font=footnotesize,labelfont=sf,textfont=sf]{subfig}
\else
  \usepackage[caption=false,font=footnotesize]{subfig}
\fi
\usepackage{url}
\begin{document}
\title{Corrigendum to "Balance of Communication and Convergence: Predefined-time Distributed Optimization Based on Zero-Gradient-Sum"}

\author{Renyongkang~Zhang,~
        Ge~Guo,~\IEEEmembership{Senior~Member,~IEEE},~
        and~Zeng-Di~Zhou% <-this % stops a space
\thanks{This paper has been accepted for publication in the IEEE Transactions on Cybernetics, doi: 10.1109/TCYB.2024.3498323.}
\thanks{This work was supported by the National Natural Science Foundation of China under Grants 62173079 and U1808205, in part by the Fundamental Research Funds for the Central Universities under Grant N2423049, and in part by the 2024 Hebei Provincial Doctoral Candidate/Postgraduate Student Innovation Ability Training Funding Project under Grant CXZZBS2024184 and CXZZSS2024178. \textit{(Corresponding author: Ge Guo.)}}
\thanks{R. Zhang and Z. D. Zhou are with the College of Information Science and Engineering, Northeastern University, Shenyang 110819, China (e-mail: zryk1998@163.com; zhouzd199912@163.com).}
%\protect\\
%E-mail: zryk1998@163.com, zhouzd199912@163.com
\thanks{G. Guo is with the State Key Laboratory of Synthetical Automation for Process Industries, Northeastern University, Shenyang 110819, China, and also with the School of Control Engineering, Northeastern University at Qinhuangdao, Qinhuangdao 066004, China (e-mail: geguo@yeah.net).}% <-this % stops an unwanted space
}


\markboth{}%
{Zhang \MakeLowercase{\textit{et al.}}: Corrigendum to "Balance of Communication and Convergence: Predefined-time Distributed Optimization Based on Zero-Gradient-Sum"}



\maketitle
\begin{abstract}
This paper proposes a distributed optimization algorithm with a convergence time that can be assigned in advance according to task requirements. To this end, a sliding manifold is introduced to achieve the sum of local gradients approaching zero, based on which a distributed protocol is derived to reach a consensus minimizing the global cost. A novel approach for convergence analysis is derived in a unified settling time framework, resulting in an algorithm that can precisely converge to the optimal solution at the prescribed time. The method is interesting as it simply requires the primal states to be shared over the network, which implies less communication requirements. The result is extended to scenarios with time-varying objective function, by introducing local gradients prediction and non-smooth consensus terms. Numerical simulations are provided to corroborate the effectiveness of the proposed algorithms.
\end{abstract}



\IEEEpeerreviewmaketitle

\section{Introduction}

\IEEEPARstart{D}{istributed} optimization has attracted considerable concern due to its widespread applications, e.g., in smart grids, sensor networks, and intelligent transportation systems\cite{ShiP2021TSMC}. The goal of distributed optimization is to cooperatively minimize a global objective formed by a summation of local cost functions using only local computation and communication. Various results on distributed optimization are available in both discrete-time (e.g., \cite{Tsitsiklis1986TAC,Nedic2009TAC,Duchi2011dual,Shi2015extra}) and continuous-time (e.g., \cite{Gharesifard2014TAC,Solmaz2015Auto,Zhang2023continuous}) domains. In particular, many continuous-time algorithms are well-known, such as (sub)gradient-based consensus algorithms \cite{Lin2016distributed,Lu2017distributed}, saddle point dynamics \cite{Gharesifard2014TAC}, proportion-integration (PI) algorithms \cite{Solmaz2015Auto,Wang2022consensus-based}, zero-gradient-sum (ZGS) schemes \cite{2012TAC-LuJ}, prediction-correction algorithms\cite{Salar2017distributed}, Newton-Rapson (NR) methods \cite{Moradian2022distributed}, etc. These algorithms can achieve asymptotic or exponential convergence, namely, the optimal solution is achieved when the time approaches infinity.

Distributed optimization algorithms ensuring convergence in a finite or fixed time are more preferable \cite{Firouzbahrami2022Cooperative,Garg2020Fixed-Time,2020Auto-LiSH,Guo2023local-minimization-free,ShiP2024TCS1}. But the upper bound of the settling time is often difficult to specify, either depending on system initialization or involving complex parameter calculation. Therefore, a great deal of attention has been paid to predefined-time optimization, which allows the maximum convergence time to be explicitly prescribed in advance. Given the challenges of the problem, some of the works give approximate solutions or alternatives that are feasible but not strictly predefined-time optimal. For instance, the method in \cite{LiKX2023predefined-time} can achieve predefined-time consensus but asymptotic optimization; the algorithms in \cite{Zhou2022distributed_INS,Zheng2023specified-time,WangYW2024TPS} can drive the system to reach a feasible domain in a given time but the optimal solution is only approachable asymptotically; the approximate solutions in \cite{Guo2020predefined-time,Liu2023multi-objective,Guo2022distributed_Auto} can converge to the neighborhood of the optimizer, but not exact optimization in a given time.

\begin{table*}[ht]
	\centering
	\caption{Comparison With Existing Predefined-Time Algorithms}
	\resizebox{\textwidth}{!}{
		\begin{tabular}{c|c|c|c|c}
			\bottomrule
			Algorithm             & \begin{tabular}[c]{@{}c@{}}Predefined-Time\\ Mechanism\end{tabular} & \begin{tabular}[c]{@{}c@{}}Predefined-Time \\ Exact Optimization\end{tabular} & Key Idea                                       & \begin{tabular}[c]{@{}c@{}}Required State of\\ Communication\end{tabular} \\ \hline
			\cite{Zhou2022distributed_INS,LiKX2023predefined-time} & Bipartite power feedback \cite{Polyakov2012Nonlinear} & No  &  Gradient-based consensus algorithm & Dimension n                                                              
			\\ \hline
			\cite{Zheng2023specified-time}   & Time-varying scaling method \cite{Pal2020design}  & No  & Prediction-correction algorithm & Dimension 2n                                                              
			\\ \hline
			\cite{Liu2023multi-objective} & Time-based generator \cite{Becerra2017predefined-time}  & No  & Saddle point dynamics & Dimension 2n
			\\ \hline
			\cite{Gong2021distributed,Cui2023Prescribed}  & \begin{tabular}[c]{@{}c@{}} \cite{Gong2021distributed}: Time-varying scaling method \cite{Wang2019Prescribed-Time} / \\ \cite{Cui2023Prescribed}: Bipartite power feedback \cite{Polyakov2012Nonlinear} \end{tabular} & Yes & Estimation + consensus + gradient dynamics & \begin{tabular}[c]{@{}c@{}}\cite{Gong2021distributed}: Dimension 2n / \\ \cite{Cui2023Prescribed}: Dimension 3n \end{tabular}
			\\ \hline
			\cite{Chen2023distributed_TCNS,De2023predefined,Ma2022distributed} & \begin{tabular}[c]{@{}c@{}}\cite{Chen2023distributed_TCNS}: Time-based generator \cite{Becerra2017predefined-time} / \\ \cite{De2023predefined}: Unilateral power feedback / \\ \cite{Ma2022distributed}: Bipartite power feedback \cite{Hu2021Fixed/Preassigned-Time} \end{tabular} & Yes & \cite{Chen2023distributed_TCNS}: ZGS / \cite{De2023predefined,Ma2022distributed}: Piecewise ZGS  & Dimension n  
			\\ \hline
			Proposed algorithm & Unilateral power feedback & Yes & Local-minimization-free ZGS & Dimension n 
			\\ \bottomrule
	\end{tabular}}
	\label{tab_Comp}
\end{table*}

%Fortunately, two types of strictly prescribed-time distributed optimization algorithms are given in some recent works. 
In some recent works, two types of strictly prescribed-time distributed optimization algorithms are given. 
The first type of methods are based on consensus and gradient dynamics with the aid of the estimation of global cost information \cite{Gong2021distributed,Cui2023Prescribed}. These algorithms can achieve the optimal solution in a preset time, but require the sharing of extra information such as estimator states or local gradients via the network, which implies poor privacy and high communication cost. The second type of methods are derived in the context of ZGS or piecewise ZGS schemes, which are computationally and communicationally more promising for predefined-time distributed optimization since only the primal states are shared among the nodes \cite{Chen2023distributed_TCNS,De2023predefined,Ma2022distributed}. It should be noted, however, that ZGS-based and piecewise ZGS-based algorithms are still limited and conservative, in that the former relies on initial conditions \cite{Chen2023distributed_TCNS} while the latter involves priority local minimization \cite{Ma2022distributed,De2023predefined}.

Motivated by the aforementioned discussions, this paper aims to develop a predefined-time distributed optimization algorithm with less demanding initial conditions and communication requirements. To this end, a local-minimization-free ZGS framework based on sliding mode is designed, consisting of a sliding manifold to guarantee that the sum of gradients converges to zero and a distributed protocol to achieve global optimal consensus. We give a uniform framework for prescribing the convergence time of both the reach phase and the sliding phase, which is very interesting due to its flexibility of regulating and fine-tuning the settling time in accordance with system performance and task objectives. This advantage is highlighted by extending the method to scenarios with time-varying objective functions, resulting in a practical algorithm running on local gradients prediction and consensus seeking. A detailed comparison with other related methods is provided in Table \ref{tab_Comp} and the main contributions are summarized as follows:

\begin{enumerate}
	\item Instead of approximate optimization \cite{Liu2023multi-objective}, predefined-time consensus-based asymptotic optimization \cite{LiKX2023predefined-time} or feasible domain convergence-based asymptotic optimization \cite{Zhou2022distributed_INS,Zheng2023specified-time}, the algorithm presented in this paper can achieve real prescribed-time distributed optimization.
	\item Superior to the methods in \cite{Gong2021distributed,Cui2023Prescribed}, in which a variety of information about the estimator and even gradients is involved, our algorithm only requires the primal states to be shared over the network, implying low communication load and strong privacy preservation.
	\item The convergence time is prescribed in a uniform framework and independent of any other parameters. Compared to the traditional predefined-time ZGS-based distributed optimization methods \cite{Chen2023distributed_TCNS,Ma2022distributed,De2023predefined}, the proposed algorithm is free from local minimization and less conservative.
\end{enumerate}

\section{Notations and Preliminaries}
%\subsection{Notations and Graph Theory}

%Denote by $ \mathbb{R} $ the set of all real numbers, and $ \mathbb{R}^+ $, $ \mathbb{R}^n $ and $ \mathbb{R}^{n\times m}  $ represent the sets of positive real numbers, real $ n $-dimension vectors, and real $ n \times m $ matrices. $ {\bf 1}_n \in \mathbb{R}^n$ and $ {\bf 0}_n \in \mathbb{R}^n$ denote an all-one column vector and an all-zero column vector, respectively. $ \boldsymbol{I}_n \in \mathbb{R}^{n\times n}$ is the identity matrix.
% Define \textcolor{blue}{$ {\rm sig}^{\alpha}\left(\cdot \right) = {\left|\cdot \right|^\alpha }\mathop{\rm sgn}\left(\cdot\right)$}, where $ \mathop{\rm sgn}\left( \cdot \right) $ is the symbolic function and $ \alpha \in \mathbb{R} $. 
Denote by $ \mathbb{R} $ the set of all real numbers, and $ \mathbb{R}^+ $, $ \mathbb{R}^n $ and $ \mathbb{R}^{n\times m}  $ represent the sets of positive real numbers, real $ n $-dimension vectors, and real $ n \times m $ matrices. $ {\bf 1}_n \in \mathbb{R}^n$ and $ {\bf 0}_n \in \mathbb{R}^n$ denote an all-one column vector and an all-zero column vector, respectively. $ \boldsymbol{I}_n \in \mathbb{R}^{n\times n}$ is the identity matrix. Define $ {\rm sig}^{\alpha}\left(\cdot \right) = {\left|\cdot \right|^\alpha }\mathop{\rm sgn}\left(\cdot\right)$, where $ \mathop{\rm sgn}\left( \cdot \right) $ is the symbolic function and $ \alpha \in \mathbb{R} $. For $ \boldsymbol{z} \in \mathbb{R}^n $ with $ \boldsymbol{z} = \left[z_1,z_2,\dots,z_n\right]^T $, $ \mathop{\rm sgn}\left( \boldsymbol{z} \right) = \left[ \mathop{\rm sgn}\left( z_1 \right),\mathop{\rm sgn}\left( z_2 \right),\dots,\mathop{\rm sgn}\left( z_n \right) \right]^T \in \mathbb{R}^n $ and $ {\rm sig}^{\alpha}\left( \boldsymbol{z} \right) = \left[ {\rm sig}^{\alpha}\left( z_1 \right),{\rm sig}^{\alpha}\left( z_2 \right),\dots, {\rm sig}^{\alpha}\left( z_n \right) \right]^T \in \mathbb{R}^n $. $ \left\|  \boldsymbol{z}  \right\| _p = \left( \left| z_1 \right|^p + \dots + \left| z_n \right|^p \right)^{1/p} $ with $p \ge 1$ represents the $ p $-norm, and the 2-norm is abbreviated as $ \left\|  {\bf z}  \right\| $. For $\boldsymbol{z} \neq {\bf 0}_n$, the norm-normalized signum function is defined as \textcolor{blue}{${\rm SGN}(\boldsymbol{z}) = \boldsymbol{z}/\|\boldsymbol{z}\| $ } and \textcolor{blue}{${\rm SIG}^\alpha(\boldsymbol{z}) = \|\boldsymbol{z}\|^\alpha {\rm SGN}(\boldsymbol{z}) $ }. Specifically, ${\rm sig}^{\alpha}\left(\boldsymbol{0} \right) =\boldsymbol{0}$, \textcolor{blue}{${\rm SGN}(\boldsymbol{{\bf 0}}) = {\bf 0}$} and \textcolor{blue}{${\rm SIG}^\alpha(\boldsymbol{0}) = \boldsymbol{0}$}. Denote by $ {\rm diag} \left( \boldsymbol{z} \right) \in \mathbb{R}^{n\times n}  $ a diagonal matrix with elements $z_i$.
 Define the eigenvalues of $ A \in \mathbb{R}^{n\times n} $ as $ \lambda _1(A) \le \cdots \le \lambda _n(A) $ in a non-decreasing order with respect to real parts. $ \nabla f\left(  \cdot  \right) $ and $ \nabla ^2 f\left(  \cdot  \right) $ respectively denote the gradient and Hessian of the twice differentiable function $ f\left(  \cdot  \right) $.

A differentiable function $ f: \mathbb{R}^n \to \mathbb{R} $ on a convex set $ {\mathcal S}\in \mathbb{R}^n $ is convex if and only if $ \forall x,z \in {\mathcal S}$, ${\left( {z - x} \right)^T}\left( {\nabla f\left( z \right) - \nabla f\left( x \right)} \right) \ge 0 $. 
For a twice continuously differentiable function $ f: \mathbb{R}^n \to \mathbb{R} $ and a position constant $\psi \in \mathbb{R}$, if $ f $ is $ \psi $-strongly convex, then the following equivalence conditions are true: $ {\left( {\nabla f\left( z \right) - \nabla f\left( x \right)} \right)^T}\left( {z - x} \right) \ge \psi {\left\| {z - x} \right\|^2}, \forall x,z \in {\mathcal S} $; $ f\left( z \right) - f\left( x \right) - \nabla f{\left( x \right)^T}\left( {z - x} \right) \ge \frac{\psi }{2}{\left\| {z - x} \right\|^2}, \forall x,z \in {\mathcal S} $; $ {\nabla ^2}f\left( x \right) \ge \psi {{\boldsymbol{I}}_n}, \forall x \in {\mathcal S} $. Besides, for any twice continuously differentiable function $ f: \mathbb{R}^n \to \mathbb{R} $ and any constant $ \varPsi >0 $, if $ \nabla f $ is $ \varPsi _i $-Lipschitz, then the following conditions are equivalent: $ {\left( {\nabla f\left( z \right) - \nabla f\left( x \right)} \right)^T}\left( {z - x} \right) \le \varPsi {\left\| {z - x} \right\|^2}, \forall x,z \in {\mathcal S} $; $ f\left( z \right) - f\left( x \right) - \nabla f{\left( x \right)^T}\left( {z - x} \right) \le \frac{\varPsi }{2}{\left\| {z - x} \right\|^2}, \forall x,z \in {\mathcal S} $; $ {\nabla ^2}f\left( x \right) \le \varPsi {{\boldsymbol{I}}_n}, \forall x \in {\mathcal S} $.

Denote by $ {\cal G} = \left( {{\cal V},{\cal E}} \right) $ a graph with the node set $ {\mathcal V} \!=\! \left\{ {1,\ldots ,N} \right\} $ and the edge set $ {\mathcal E} \!\subseteq \!{\mathcal V} \!\times\! {\mathcal V} $. Define $ {\mathcal A} \!=\!\left[a_{ij}\right] \!\in\! \mathbb{R}^{N\!\times\! N} $ as the weighted adjacency matrix with $ {a_{ii}} \!=\! 0 $, $ {a_{ij}} \!>\! 0 $ if $ \left( {j,i} \right) \!\in\! {\mathcal E} $ and $ {a_{ij}} \!=\! 0 $ otherwise. $ \mathcal{L} \!=\! \left[ {{l_{ij}}} \right] \!\in\! \mathbb{R}^{N\times N} $ represents the Laplacian matrix of $ {\mathcal G} $, where $ {l_{ij}} \!=\! - {a_{ij}} $ for $ i \!\ne\! j $ and $ {l_{ii}} \!=\! \sum_{j = 1}^N {{a_{ij}}} $. For any undirected and connected graph $ {\mathcal G} $, the following conclusion holds: 1) $ {\mathcal L} $ is positive semidefinite and satisfies $ {\lambda _1}\left( {\mathcal L} \right) = 0 < {\lambda _2}\left( {\mathcal L} \right) \le  \cdots  \le {\lambda _N}\left( {\mathcal L} \right) $, 2) for $ x=[x_i]^T \in \mathbb{R}^n $ and $ y \in \mathbb{R}^n $, it holds that $ x^T \mathcal{L} x = \frac{1}{2} \sum_{i,j}a_{ij}\left( x_i - x_j\right)^2$ and $ \sum_{i,j}a_{ij} x_i \left( y_i - y_j\right) = \frac{1}{2} \sum_{i,j}a_{ij}\left( x_i - x_j\right)\left( y_i - y_j\right) $, 3) if $ {{\bf 1}_n}^T x =0 $, $ x^T \mathcal{L} x \ge {\lambda _2}\left( \mathcal{L} \right) x^Tx $.

%Let $ {\cal G} = \left( {{\cal V},{\cal E},{\cal A}} \right) $ denote a graph with the node set $ {\cal V} = \left\{ {1,2, \ldots ,N} \right\} $, the edge set $ {\cal E} \subseteq {\cal V} \times {\cal V} $ and the weighted adjacency matrix $ {\cal A} = \left[ {{a_{ij}}} \right] \in \mathbb{R}^{N\times N} $, where $ {a_{ii}} = 0 $, $ {a_{ij}} =1 $ if $ \left( {j,i} \right) \in {\cal E} $ and $ {a_{ij}} = 0 $ otherwise. $ {\cal L} = \left[ {{l_{ij}}} \right] \in \mathbb{R}^{N\times N} $ denotes the Laplacian matrix of $ {\cal G} $, where $ {l_{ij}} = - {a_{ij}} $ for $ i \ne j $ and $ {l_{ii}} = \sum_{j = 1}^N {{a_{ij}}} $. The neighbor set of node $ i $ is $ {{\cal N}_i} = \left\{ {j \in {\cal V}\left| {\left( {j,i} \right) \in {\cal E} } \right.} \right\} $. The path of $ {\cal G} $ is an edge sequence $ \left( {{k_1},{k_2}} \right),\left( {{k_2},{k_3}} \right), \ldots ,\left( {{k_{p - 1}},{k_p}} \right) $ for $ {k_1}, \ldots ,{k_p} \in {\cal V} $. A graph is connected if there is a path between any two distinct nodes. Generally, if $ \forall i,j \in {\cal V} $, $ {a_{ij}} = {a_{ji}} $, the $ {\cal G} $ is said to be undirected graph and $ 1_N^T{\cal L} = 0 $. In connected undirected graph $ {\cal G} $, its Laplacian matrix $ {\cal L} $ is positive semidefinite, and eigenvalues of $ {\cal L} $ are denoted by $ {\lambda _1}\left( {\cal L} \right) < {\lambda _2}\left( {\cal L} \right) \le  \cdots  \le {\lambda _N}\left( {\cal L} \right) $, where $ {\lambda _1}\left( {\cal L} \right) $ is zero with associated eigenvector $ 1_N $. 

%\subsection{Definitions and Lemmas}
\begin{lemma}[\cite{Zuo2014IJC}]
	\label{lemma_Nn}
	For $ {z_i} \in \mathbb{R} \ge 0 $ and $ p \in \mathbb{R}$, if $ 0 < p \le 1 $, then $ \sum\nolimits_{i = 1}^n {z_i^p}  \ge {\left( {\sum\nolimits_{i = 1}^n {{z_i}} } \right)^p} $, and if $ p > 1 $, then $ \sum\nolimits_{i = 1}^n {z_i^p}  \ge {n^{1 - p}}{\left( {\sum\nolimits_{i = 1}^n {{z_i}} } \right)^p} $. 
\end{lemma}
\begin{lemma}[\cite{Sanchez2020predefined-time_IJC}]
	\label{lemma_PdT-time0}
	Consider the system $ \dot x\left( t \right) \!=\! f\left( {x\left( t \right)} \right) $,
	where $ x \!\in\! \mathbb{R}^n $ and $ f: \mathbb{R}^n \to \mathbb{R}^n $. If there exists a Lyapunov function $ {V}\left( x \right)$ such that $$ \dot {V}\left(x \right) \le  - \frac{\alpha^{\frac{\beta q-1}{p}} \varGamma\left(\frac{1-\beta q}{p}\right)}{ p T_m} \exp(\alpha {V}^p ){V}^{\beta q}, \forall x \!\in\! \mathbb{R}^n / \{\textbf{0}\}, $$ where $T_m >0$, $ \alpha >0 $, $ p>0 $, $\beta >0$, $q>0$, $\beta q<1$, and $\varGamma\left(\tau \right) = \int_{0}^{\infty}t^{\tau - 1}\exp\left(-t\right)dt$, then the equilibrium of system is globally predefined-time stable and the settling time can be bounded by $ T_m $.
\end{lemma}
\begin{lemma}
	\label{lemma_PdT-time}
	For the system $ \dot x\left( t \right) \!=\! f\left( {x\left( t \right)} \right) $ with $ x \!\in\! \mathbb{R}^n $ and $ f: \mathbb{R}^n \to \mathbb{R}^n $, if there exists a Lyapunov function $ {V}\left( x \right)$ such that $ \dot {V}\left(x \right) \le  - \frac{1}{\alpha p T_m} \exp(\alpha {V}^p ){V}^{1-p}, \forall x \!\in\! \mathbb{R}^n / \{\textbf{0}\} $, where $T_m >0$, $ \alpha >0 $ and $ 0 < p < 1 $, then the equilibrium of system is globally predefined-time stable and the settling time is bounded by $ T_m $.
\end{lemma}

\begin{IEEEproof}
	By Lemma \ref{lemma_PdT-time0}, let $\beta q=1-p$, $\beta =1$, and $ 0<p<1$,  the conclusion can be drawn.
\end{IEEEproof}

\begin{remark}
	Unlike the bipartite power feedback in \cite{Cui2023Prescribed}, which contains both powers greater than zero but less than one and powers greater than one, Lemma \ref{lemma_PdT-time} achieves predefined-time stability by only using unilateral power feedback of powers greater than zero but less than one, reducing computational complexity. A distributed optimization algorithm with predefined-time convergence based on Lemma \ref{lemma_PdT-time} will be presented later in the article.
\end{remark}

\section{Problem Formulation}
Consider a system with $N$ agents interacting over a communication topology $ \cal G  $. Each agent $ i \in {\cal V} $ is endowed with a local cost function $ {f_i}\left( {{z}} \right): \mathbb{R}^n \to \mathbb{R} $ or a time-varying cost function $ {f_i}\left( {{z}},t \right): \mathbb{R}^n \times \mathbb{R}^+ \to \mathbb{R} $, which is only accessible to itself. The aim of this paper is to design a distributed algorithm for each agent to solve the global optimization problem
$ \min \sum\nolimits_{i = 1}^N {{f_i}\left( {{z}} \right)}$ or the global time-varying optimization $ \min\sum\nolimits_{i = 1}^N {{f_i}\left( {{z},t} \right)} $ using only local information and interaction with neighbors.
%\begin{align}
%	\begin{array}{c}
	%		\min\limits_{{z}} \sum\limits_{i = 1}^N {{f_i}\left( {{z}} \right)} \\
	%	\end{array}
%\end{align}
%or the global time-varying optimization
% \begin{align}
	% 	\begin{array}{c}
		% 		\min\limits_{{z}} \sum\limits_{i = 1}^N {{f_i}\left( {{z},t} \right)} \\
		% 	\end{array}
	% \end{align}
Since the decision variables of each agent are the same, the above problem can be formulated as
\begin{align}
	\label{optimization_problem_Time-Invariant}
	\begin{array}{c}
		\min\limits_{x_i} \sum\limits_{i = 1}^N {{f_i}\left( {{x_i}} \right)}~~~~{\text{s.t.}}~~{x_i} = {x_j}~~\forall i,j \in {\mathcal V}
	\end{array}
\end{align}
or 
\begin{align}
	\label{optimization_problem_Time-varying}
	\begin{array}{c}
		\min\limits_{x_i} \sum\limits_{i = 1}^N {{f_i}\left( {{x_i},t} \right)}~~~~{\text{s.t.}}~~{x_i} = {x_j}~~\forall i,j \in {\mathcal V}
	\end{array}
\end{align}
%${\bf x}=[x_1^T,\cdots,x_N^T]^T \in \mathbb{R}^{nN} $, and
where $x_i \in \mathbb{R}^n $ is the state of the $i$th agent. Suppose that the optimization problem \eqref{optimization_problem_Time-Invariant} and \eqref{optimization_problem_Time-varying} are feasible and denote by $x^*$ the optimal solution. Note that $x^*$ is a static value for the problem \eqref{optimization_problem_Time-Invariant} while $x^*$ is a continuous trajectory that changes over time for \eqref{optimization_problem_Time-varying}. For the sake of convenience, the gradients and Hessian of $ f_i( \cdot) $ with respect to $x_i$ will be denoted as $ \nabla f_i (\cdot ) $ and $ \nabla ^2 f_i (\cdot ) $, respectively, throughout the following text.

Next several assumptions and definitions are required in following analysis.
\begin{assumption}
	\label{assumption_G}
	$ \mathcal G $ is undirected and connected.
\end{assumption}

\begin{assumption}
	\label{assumption_cost}
	Each local cost function $ f_i $ is twice continuously differentiable, and $ \psi_i $-strongly convex ($ \psi_i > 0$) with respect to $x_i$. In addition, $ \nabla f_i $ is $ \varPsi _i $-Lipschitz ($ \varPsi _i >0$) continuous with respect to $x_i$. 
\end{assumption}

\begin{remark}
	Assumption \ref{assumption_G} and \ref{assumption_cost} both are common standard assumptions and often used in the existing literature. The existence and uniqueness of the solution of the optimization problem is guaranteed by Assumption \ref{assumption_cost}, and the optimal solution $ x^* $ satisfying $ \sum_{i = 1}^N{ \nabla f_i\left(  x^*  \right)} = 0 $. %$ \psi _i \boldsymbol{I}_n \le \nabla ^2 f_i\left( \cdot \right) \le \varPsi _i \boldsymbol{I}_n$
\end{remark}
%\begin{remark}
%	Assumptions \ref{assumption_G} and \ref{assumption_cost} are common standard assumptions in distributed optimization. From \cite{Olfati-Saber2004TAC}, if Assumption \ref{assumption_G} holds, one has $ \xi^T\! (\mathcal{L} \otimes \boldsymbol{I}_n) \xi \!=\! \frac{1}{2} \sum_{i,j}a_{ij}\left\| \xi_i \!-\! \xi_j\right\|^2$ and $ \xi^T \! (\mathcal{L}\otimes \boldsymbol{I}_n) \xi \!\ge\! {\lambda _2}\left( \mathcal{L} \right) \xi^T \! \xi $ if $ {\bf 1}_{nN}^T\! \xi \!=\!{\bf 0} $, where $ \xi=[\xi_1^T,\xi_2^T,\dots,\xi_N^T]^T \!\in\! \mathbb{R}^{nN} $ and $ \xi_i \!\in\! \mathbb{R}^{n}, i \in {\mathcal V}$. Assumption \ref{assumption_cost} implies the positive definiteness of Hessian $ \nabla ^2 f_i\left(  x_i  \right) $ and the uniqueness of the global optimizer $ x^* $, which satisfies $ \sum_{i = 1}^N{ \nabla f_i\left(  x^*  \right)} \!=\! {\bf 0} $.
%\end{remark}

\begin{definition}[Predefined-time optimization]
	For the problem \eqref{optimization_problem_Time-Invariant} or \eqref{optimization_problem_Time-varying}, it is said to achieve predefined-time optimization if there exists a continuous-time algorithm such that the optimal solution $ x^* $ is sought at a predefined time $T_m$, i.e., $ \lim\limits_{t \to T_m} \left\| x_i(t) - x^* \right\|  = 0  $ and $\left\| x_i(t) - x^* \right\| = 0 $ for $ \forall t > T_m$.
\end{definition}
\begin{definition}[Predefined-time approximate optimization,\cite{Liu2023multi-objective}]
	For the problem \eqref{optimization_problem_Time-Invariant} or \eqref{optimization_problem_Time-varying}, it is said to achieve predefined-time approximate optimization if there exists a continuous-time algorithm such that
	\begin{itemize}
		\item $\left\| x_i(t) - x^* \right\| \!\le\! \varepsilon $ for $ \forall t \!>\! T_m$ and $ \lim_{t \to T_m^+} \left\| x_i(t) - x^* \right\| \!\le\! \varepsilon $.
		\item  $ \lim_{t \to \infty} \left\| x_i(t) - x^* \right\| = 0 $.
	\end{itemize}
\end{definition}

\section{Main Results}

\subsection{Predefined-time Distributed Optimization Algorithms}
%\subsubsection{Predefined-time Algorithms}
In this section, a novel predefined-time distributed optimization algorithm is first provided to solve the time-invariant problem stated in \eqref{optimization_problem_Time-Invariant}. For each agent $i \in \mathcal{V}$ and positive constants $ \eta, p , T_m, c $,  consider
\begin{subequations}
	\label{Predefined_Algorithm}
	\begin{align}
		\label{PdT_Algorithm_dynamic}
		\dot{x}_i(t) = & - \big({ {\nabla ^2}{f_i}\left( {{x_i}} \right)}\big)^{-1} \Big(  \frac{1}{2p\eta T_m}\exp(\left\| s_i \right\|^{2p} ){\rm sig}^{1-2p}\left( s_i \right)  \notag \\
		& + \frac{2c}{p(1-\eta)T_m} \sum\limits_{j = 1}^N  \exp((a_{ij}\left\| {{x_i} - {x_j}} \right\|^2 )^p) \times \notag \\
		& a_{ij}^{1-p}{\rm sig}^{1-2p}\left( {x_i} - {x_j} \right)  \Big), \\
		\label{PdT_Algorithm_sliding mode}
		{s_i}(t) =&~ \nabla \!{f_i}\left(\! {{x_i}} \!\right) \!+\!\! \int_0^t \!\! \frac{2c}{p(1\!-\! \eta)T_m} \Big(\! \sum\limits_{j = 1}^N \! \exp((a_{ij}\left\| \! {{x_i}\! -\! {x_j}} \right\|^2 )^p) \times \notag \\
		& a_{ij}^{1-p}{\rm sig}^{1-2p}\left( {x_i} - {x_j} \right)  \Big) d\tau 
	\end{align}
\end{subequations}
\eqref{PdT_Algorithm_dynamic} is a continuous-time distributed update law to solve optimization problem, and \eqref{PdT_Algorithm_sliding mode} is a sliding manifold to achieve zero-gradient-sum. The convergence analysis of algorithm \eqref{Predefined_Algorithm} is given as the following theorem.

\begin{theorem}
	\label{theorom_PT_Algorithm}
	Suppose that Assumptions \ref{assumption_G} and \ref{assumption_cost} hold. Let $ \bar{\varPsi} = \max_{ i \in {\mathcal V} } \{\varPsi _i\} $. If $ 0< \eta <1 $, $ 0 < p < {1}/{2}$, $ T_m >0 $, and $ c \ge {\bar{\varPsi}}/\left({4 \lambda _2 \left( \mathcal{L} \right)}\right) $, then, distributed algorithm \eqref{Predefined_Algorithm} drives the states $x_i(t)$ of all agents to a optimal solution $x^*$ of \eqref{optimization_problem_Time-Invariant} within a settling time that is bounded by $ T_m $. 
\end{theorem}

\begin{IEEEproof}
%\emph{Proof:}
The proof is divided into two steps. Step 1: To prove that the sum of local gradients goes to zero in $\eta T_m$; Step 2: To show that the states $ x_i(t) $ converges to the optimal solution within $T_m$.

\textbf{Step 1:} For $ 0< t \le \eta T_m$, the time derivative of \eqref{PdT_Algorithm_sliding mode}
satisfies
\begin{align}
	\label{PdT_Algorithm_sliding mode_dynamic}
	\dot s _i (t) = - \frac{1}{2p\eta T_m}\exp(\left\| s_i \right\|^{2p} ) {\rm sig}^{1-2p}\left( s_i \right)
\end{align}
Consider the Lyapunov function candidate ${V} _{s,i} \left(s_i(t)\right) = \frac{1}{2}{\left\|  s_{i}  \right\|^2}$. The time derivative of $ {V}_{s,i}\left(t\right) $ along \eqref{PdT_Algorithm_sliding mode_dynamic} is given by $\dot{V} _{s,i} = - \frac{1}{2p\eta T_m}\exp(\left\| s_i \right\|^{2p} ){\left\| s_i \right\|^{2-2p}_{2-2p}} $.
Since $0<p<1/2$, one has form Lemma \ref{lemma_Nn} that
\begin{align}
	\dot{V} _{s,i} \le& - \frac{1}{2p\eta T_m}\exp(\left\| s_i \right\|^{2p} ){({\left\| s_i \right\|^2})^{1-p}} \notag \\
	=& - \frac{1}{{2^p}p\eta T_m}\exp(2^p V_{s,i}^{p} ){V_{s,i}^{1-p}}
\end{align}
By Lemma \ref{lemma_PdT-time}, it can be obtained that $ \!s_i(t) , i \!\in\! \mathcal{V}\!  $ approach zero in a specified time $ \eta T_m $. Since Assumption \ref{assumption_G} holds, it follows that
\begin{align}
	\sum\limits_{i = 1}^N \sum\limits_{j = 1}^N { \exp((a_{ij}\left\| {{x_i} - {x_j}} \right\|^2 )^p) a_{ij}^{1-p}{\rm sig}^{1-2p}\left( {x_i} - {x_j} \right) }= \textbf{0} 
\end{align}
From \eqref{PdT_Algorithm_sliding mode}, we have
\begin{align}
	\sum\nolimits_{i = 1}^N {\nabla \!{f_i}\!\left( {{x_i(t)}} \right)} = \sum\nolimits_{i = 1}^N {{s_i(t)}}
\end{align}
Combining with $ \lim_{t \to \eta T_m} s_i \left( t \right) = {\bf 0} $, we can get $ \sum\nolimits_{i = 1}^N \!{\nabla \!{f_i}\!\left( {{x_i(t)}} \right)}  \!=\! \sum\nolimits_{i = 1}^N \!{{s_i(t)}} \! =\!  {\bf 0} $ at $ t \!=\! \eta T_m $.

\textbf{Step 2:}
For $ t > \eta T_m$, it follows from $ s_i ={\bf 0} $ that the dynamics satisfies 
\begin{align}
	\label{PdT_Algorithm_closed_loop_dynamic}
	\dot{x}_i (t) =& -\! \big({\! {\nabla ^2}\!{f_i}\left(\! {{x_i}}\! \right)}\big)^{-1} \!\Big(\! \frac{2c}{p(1\!-\!\eta)T_m}\!\! \sum\limits_{j = 1}^N \! \exp(\!(a_{ij}\!\left\| {{x_i} \!-\! {x_j}} \!\right\|^2 )^p)\!\times \notag \\
	& a_{ij}^{1-p}{\rm sig}^{1-2p}\left( {x_i} - {x_j} \right)  \Big)
\end{align}
Substituting \eqref{PdT_Algorithm_closed_loop_dynamic}, it can be derived that ${{d \sum\nolimits_{i = 1}^N {\nabla {f_i}\left( {{x_i}(t)} \right)}}}/{{dt}} = \sum\nolimits_{i = 1}^N {\nabla^2 {f_i}\left( {{x_i}(t)} \right)} \dot{x}_i = {\bf 0} $. Due to $ \sum\nolimits_{i = 1}^N {\nabla {f_i}\left( {{x_i}(t)} \right)}  = {\bf 0} $ at $ t = \eta T_m $, it holds that $ \sum\nolimits_{i = 1}^N {\nabla {f_i}\left( {{x_i}(t)} \right)}  = {\bf 0}, \forall t \ge \eta T_m $, which indicates that the zero-gradient-sum manifold is guaranteed.

Choose the Lyapunov function for system \eqref{PdT_Algorithm_closed_loop_dynamic}
\begin{align}
	V \left( x_i(t) \right) = \sum\limits_{i = 1}^N \!\! \left( {{f_i}\left( {{x^*}} \right) \!-\! {f_i}\left( {{x_i}} \right) \!-\! \nabla f_i^T\left( {{x_i}} \right)\left( {{x^*} - {x_i}} \right)} \right)
\end{align}
It is clear from Assumption \ref{assumption_cost} that $ V \!\! \ge \!\! \sum\nolimits_{i = 1}^N  \!\! \frac{\psi _i}{2}\!{\left\| \! {x_i} \! - \! {x^*} \! \right\|^2} $.
% and $ V \! = \! 0 \! \Leftrightarrow \! {x_i} \! = \! {x^*} $.
Denote $ {\bf \tilde x} = \left[\tilde x_1^T,\ldots, \tilde x_N^T \right]^T \in \mathbb{R}^{nN} $, where ${\tilde x_i} = x_i - {\bar x} $ with ${\bar x} = \frac{1}{N}\sum_{k=1}^{N}x_k $.
Since $ F\left(x^*\right) \le F\left(\bar{x}\right) $ and $\sum\nolimits_{i = 1}^N {\nabla {f_i}\left( {{x_i}(t)} \right)}  = {\bf 0}$, it is obvious that
\begin{align}
	\label{inequalityF-V}
	&\sum\limits_{i = 1}^N \!\! \left( {{f_i}\left( {\bar{x}} \right) \!-\! {f_i}\left( {{x_i}} \right) \!-\! \nabla f_i^T\left( {{x_i}} \right)\left( {\bar{x} - {x_i}} \right)} \right) - V(t) \notag \\
	&= \sum_{i=1}^{N}{ \big( {f_i}\left( {\bar{x}} \right) -{f_i}\left( {x^*} \right) \big) } - \big( \sum_{i=1}^{N}{\nabla f_i \left( {{x_i}} \right)} \big)^T \left( {\bar{x} - {x_i}} \right) \notag \\
	&= \sum_{i=1}^{N}{ \big( {f_i}\left( {\bar{x}} \right) -{f_i}\left( {x^*} \right) \big) } \ge {\bf 0}
\end{align}
Notice that the following inequalities holds
\begin{align}
	\label{inequalityF}
	&\sum\limits_{i = 1}^N \!\! \left( {{f_i}\left( {\bar{x}} \right) \!-\! {f_i}\left( {{x_i}} \right) \!-\! \nabla f_i^T\left( {{x_i}} \right)\left( {\bar{x} - {x_i}} \right)} \right) \notag \\
	&\le \sum_{i=1}^{N}{ \frac{\varPsi _i}{2} \| \tilde x _i \|^2 } \le \frac{\bar{\varPsi}}{2} \| {\bf \tilde x} \|^2
\end{align}
From \eqref{inequalityF-V} and \eqref{inequalityF}, it follows that
\begin{align}
	\label{inequalityV_X}
	V(t) \le \frac{\bar{\varPsi}}{2} \| {\bf \tilde x} \|^2
\end{align}

Differentiating  $ V \left(t\right) $ along \eqref{PdT_Algorithm_closed_loop_dynamic} yields
\begin{align}
	{\dot V } =& \sum_{i=1}^{N} \left( x_i - x^* \right)^T { {\nabla ^2}{f_i}\left( {{x_i}} \right)} {\dot{x}_i} \notag \\
	=&- \frac{2c}{p(1-\eta)T_m} \sum\limits_{i = 1}^N \sum\limits_{j = 1}^N  \Big( \exp\big((a_{ij}\left\| {{x_i} - {x_j}} \right\|^2 )^p\big)  \times \notag \\
	&  a_{ij}^{1-p} \left( x_i -x^* \right)^T {\rm sig}^{1-2p}\left( {x_i} - {x_j} \right)  \Big)
\end{align}
Given that Assumption 1 holds, one has
\begin{align}
&\sum\limits_{i = 1}^N \sum\limits_{j = 1}^N  \Big(\! \exp\big((a_{ij}\left\| {{x_i} - {x_j}} \right\|^2 )^p\big)  a_{ij}^{1-p} x_i ^T {\rm sig}^{1-2p}\left( {x_i} - {x_j} \right)  \!\Big)  \notag \\
&=\frac{1}{2} \sum\limits_{i = 1}^N \sum\limits_{j = 1}^N  \Big( \exp\big((a_{ij}\left\| {{x_i} - {x_j}} \right\|^2 )^p\big)  a_{ij}^{1-p} \left( x_i -x_j \right)^T \notag \\
&\times {\rm sig}^{1-2p}\left( {x_i} - {x_j} \right)  \Big)
\end{align}
and
\begin{align}
&\sum\limits_{i = 1}^N \sum\limits_{j = 1}^N \! \Big(\! \exp\big(\!(a_{ij}\left\| {{x_i} - {x_j}} \right\|^2 )^p\!\big) a_{ij}^{1-p} (x^*) ^T {\rm sig}^{1-2p}\left( {x_i} - {x_j} \!\right)  \! \!\Big)  \notag \\
&=0
\end{align}
It follows that
\begin{align}
	{\dot V } =&- \frac{c}{p(1-\eta)T_m} \sum\limits_{i = 1}^N \sum\limits_{j = 1}^N  \Big( \exp\big((a_{ij}\left\| {{x_i} - {x_j}} \right\|^2 )^p\big)  \times \notag \\
	&  a_{ij}^{1-p} \left( x_i -x_j \right)^T {\rm sig}^{1-2p}\left( {x_i} - {x_j} \right)  \Big)  \notag \\
	=&- \frac{c}{p(1-\eta)T_m} \sum\limits_{i = 1}^N \sum\limits_{j = 1}^N  \Big( \exp\big((a_{ij}\left\| {{x_i} - {x_j}} \right\|^2 )^p\big)  \times \notag \\
	&  a_{ij}^{1-p} \left\| x_i -x_j \right\|^{2-2p}_{2-2p} \Big)  \notag \\
	\le&- \frac{c}{p(1-\eta)T_m} \sum\limits_{i = 1}^N \sum\limits_{j = 1}^N  \Big( \exp\big((a_{ij}\left\| {{x_i} - {x_j}} \right\|^2 )^p\big)  \times \notag \\
	&  { \big( {a_{ij}} \left\| {{x_i} - {x_j}} \right\|^2 \big)^{1-p} }  \Big)
\end{align}
Note that the function $f\left(z\right) = \exp(z^{2p})z^{2(1-p)} $ is convex. By applying Jensen's Inequality, it can be deduced that
\begin{align}
&\frac{1}{N\!(\!N\!-\!1\!)} \!\! \sum\limits_{i = 1}^N \!\sum\limits_{j = 1}^N  \!\Big(\! \exp\big(\!(\!a_{ij}\left\| {{x_i} - {x_j}} \right\|^2 \!)^p\!\big)  \! { \big(\! {a_{ij}} \left\| {{x_i} - {x_j}} \right\|^2 \!\big)^{1-p} }  \!\Big) \notag \\
& \ge \exp\!\big( \!( \!\frac{\sum\limits_{i = 1}^N \!\sum\limits_{j = 1}^N a_{ij}\left\| {{x_i} - {x_j}} \right\|^2}{N(N-1)} )^p\big)\! { \big( \!\frac{ \sum\limits_{i = 1}^N\!\sum\limits_{j = 1}^N {a_{ij}} \left\| {{x_i} - {x_j}} \right\|^2}{N(N-1)} \!\big)^{1-p} } 
\end{align}
Then, 
\begin{align}
	{\dot V } \le&- \frac{cN(N-1)}{p(1-\eta)T_m} \Big( \exp\big( ( \frac{\sum\limits_{i = 1}^N \sum\limits_{j = 1}^N a_{ij}\left\| {{x_i} - {x_j}} \right\|^2}{N(N-1)} )^p\big)  \times \notag \\
	& { \big( \frac{ \sum\limits_{i = 1}^N \sum\limits_{j = 1}^N {a_{ij}} \left\| {{x_i} - {x_j}} \right\|^2}{N(N-1)} \big)^{1-p} }  \Big) \notag \\
	=&- \frac{cN(N-1)}{p(1-\eta)T_m} \Big( \exp\big( ( \frac{\sum\limits_{i = 1}^N \sum\limits_{j = 1}^N a_{ij}\left\| {{\tilde x_i} - {\tilde x_j}} \right\|^2}{N(N-1)} )^p\big)  \times \notag \\
	& { \big( \frac{ \sum\limits_{i = 1}^N \sum\limits_{j = 1}^N {a_{ij}} \left\| {{\tilde x_i} - {\tilde x_j}} \right\|^2}{N(N-1)} \big)^{1-p} }  \Big) 
\end{align}
Since the graph is undirected and connected, ${\bf 1}^T {\bf \tilde x} =0 $ and $ {\lambda _2 \left( \mathcal{L} \right)} = \lambda _2 \left( \mathcal{L} \otimes I_n\right) $, we have 
\begin{align}
	 \sum\limits_{i = 1}^N \sum\limits_{j = 1}^N a_{ij}\left\| {{x_i} - {x_j}} \right\|^2 =  2{\bf \tilde x}^T \left( \mathcal{L} \otimes I_n\right) {\bf \tilde x} \ge 2{\lambda _2 \left( \mathcal{L} \right)} \| {\bf \tilde x} \|^2
\end{align}
Considering that $f\left(z\right) = \exp(z^{2p})z^{2(1-p)} $ is monotonically increasing on $ \left[ 0, \infty \right) $, it can be concluded that
\begin{align}
	{\dot V } \le&- \frac{cN(N-1)}{p(1-\eta)T_m} \Big( \exp\big( ( \frac{2{\bf \tilde x}^T \left( \mathcal{L} \otimes I_n\right) {\bf \tilde x} }{N(N-1)} )^p\big)  \times \notag \\
	&  { \big( \frac{ 2{\bf \tilde x}^T \left( \mathcal{L} \otimes I_n\right) {\bf \tilde x} }{N(N-1)} \big)^{1-p} }  \Big) \notag \\
	\le&- \frac{cN(N-1)}{p(1-\eta)T_m} \Big( \exp\big( ( \frac{2{\lambda _2 \left( \mathcal{L} \right)} \| {\bf \tilde x} \|^2 }{N(N-1)} )^p\big)  \times \notag \\
	\label{inequality_dotV_X}
	&  { \big( \frac{ 2{\lambda _2 \left( \mathcal{L} \right)} \| {\bf \tilde x} \|^2 }{N(N-1)} \big)^{1-p} }  \Big)
\end{align}
From \eqref{inequalityV_X} and \eqref{inequality_dotV_X}, one has 
\begin{align}
	{\dot V } \le- \!\frac{cN(N\!-\!1)}{p(1\!-\!\eta)T_m} \big(\! \frac{ 4{\lambda _2 \left( \!\mathcal{L} \!\right)} }{N(N \!-\!1){\bar{\varPsi}}} \big)^{1\!-\!p} \! \exp\!\big( \!(\! \frac{4{\lambda _2 \left( \mathcal{L} \right)} V }{N(N \!-\!1){\bar{\varPsi}}} \!)^p\big)  { V^{1-p} }
\end{align}
Due to $ c \ge {\bar{\varPsi}}/{(4 \lambda _2 \left( \mathcal{L} \right))} $, it can be obtained that 
\begin{align}
	{\dot V } \le- \!\frac{1}{p(1\!-\!\eta)T_m} \!\big(\! \frac{ 4{\lambda _2 \left( \!\mathcal{L}\! \right)} }{N(N\!-\!1){\bar{\varPsi}}} \!\big)^{\!-p} \!\exp \! \big( \! (\!\frac{4{\lambda _2 \left(\mathcal{L} \right)} V }{N(N\!-\!1){\bar{\varPsi}}} \!)^p\big)  { V^{1\!-\!p} }
\end{align}
Invoking Lemma \ref{lemma_PdT-time}, we have $ V =0$ with a predefined time $(1-\eta)T_m$. With the help of optimality conditions, it is clear that $ \lim_{t \to T_m} x_i(t) \to x^*, ~ i \in {\mathcal V}$.
\end{IEEEproof}
%$\hfill\blacksquare$

\begin{remark}
The proposed algorithm is inspired by the ZGS framework \cite{Guo2023local-minimization-free}  and predefined-time stability \cite{De2023predefined}. Compared with \cite{Guo2023local-minimization-free}, the convergence time in this paper can be arbitrarily specified by users and does not depend on system parameters, and is established under the unified time by introducing a tuning parameter $\eta$. In contrast to \cite{De2023predefined}, our algorithm is less conservative, not only circumventing the requirement for local minimization, thereby avoiding waste of computational time, but also offering a broader range of parameters selection ($ c \ge {\bar{\varPsi}}/\left({4 \lambda _2 \left( \mathcal{L} \right)}\right) $ in this paper but $ k \ge {\bar{\varPsi}}/\left({2 \lambda _2 \left( \mathcal{L} \right)}\right) $ in \cite{De2023predefined} ). It is also worth mentioning that the proposed method eliminates the singularity of the algorithm \cite{De2023predefined} as the gradient approaches zero.
\end{remark}

\begin{remark}[Convergence of \eqref{Predefined_Algorithm} Over Dynamically Changing Topologies] The convergence result in Theorem \ref{theorom_PT_Algorithm} is established for fixed undirected connected graphs. Since the proof of Theorem \ref{theorom_PT_Algorithm} relies on a Lyapunov function with no dependency on the communication topology, we can readily extend the convergence result to dynamically changing topologies, following the line of proof Theorem 1. In other words, for a time-varying graph $\mathcal{G}(t)$ which is undirected and connected at all times and whose adjacency matrix is uniformly bounded and piecewise constant, the algorithm \eqref{Predefined_Algorithm} is still effective as long as $ c \ge {\bar{\varPsi}}/({4  {\lambda_2}_\text{min}  \left( \mathcal{L}(t) \right) }) $ with ${\lambda_2}_\text{min}  \left( \mathcal{L}(t) \right) = \min {\lambda_2}  \left( \mathcal{L}(t) \right)$.
\end{remark}

\begin{remark}[Free-Will Arbitrary Time Convergence] The specified time distributed optimization framework proposed in this paper relies on the predefined-time stability in \cite{Sanchez2020predefined-time_IJC}. Different convergence results can be obtained easily by using various stability lemmas. For example, based on the result in \cite{Pal2020design,Pal2020Free-will}, the algorithm with arbitrary time convergence can be obtained, as shown below:
	\begin{subequations}
		\label{FreeWill-T_Algorithm}
		\begin{align}
			\dot{x}_i(t) =& - \!\big(\!{ {\nabla ^2}{f_i}\left(\! {{x_i}}\! \right)}\big)^{-1} \!\bigg( \! \frac{k}{\eta T_m-t} \!\Big( \!I_n - \exp\!\big(\!-{\rm diag}(s_i) \!\big) \!\Big) \! \textbf{1}_n \notag \\
			&+ \frac{c}{(1-\eta)T_m -t } \sum\limits_{j = 1}^N {  {a_{ij}}\left( {{\chi_i} - {\chi_j}} \right) } \! \bigg), \notag \\
			{\chi_i}(t) =&~ \bigg( I_n - \exp \Big(-{\rm diag}\big( \sum\limits_{j = 1}^N {  {a_{ij}}\left( {{x_i} - {x_j}} \right) } \big) \Big) \bigg) \textbf{1}_n, \notag \\
			{s_i}(t) =&~ \nabla {f_i}\left( {{x_i}} \right) \!+\! \int_0^t \! \frac{c}{(1-\eta)T_m \!-\!\tau } \sum\limits_{j = 1}^N {  {a_{ij}}\left( {{\chi_i} - {\chi_j}} \right) }  d\tau 
		\end{align}
	\end{subequations}
	where $ 0< \eta <1 $, $ k \ge 1$, $ T_m >0 $, and $ c \ge {\bar{\varPsi}} / ({ \lambda ^2 _2 \left( \mathcal{L} \right)}) $. 
\end{remark}

\begin{remark}[Robustness to Disturbance] The proposed method based on sliding mode can reject disturbance after proper extension. Suppose that $ d_i \left( t \right) \in \mathbb{R}^n$ is unknown disturbance and $ \left\| d_i(t) \right\|_1 \le \mathcal{D}$. Introducing a non-smooth term $ k {\rm sgn}\left(s_i\right) $, it can be obtained the following algorithm insensitive to disturbance with predefined-time optimization.
	\begin{subequations}
		\label{PdT_Algorithm_Disturbance}
		\begin{align}
			\dot{x}_i(t) =& - \big({ {\nabla ^2}{f_i}\left( {{x_i}} \right)}\big)^{-1} \Big(  \frac{1}{2p\eta T_m}\exp(\left\| s_i \right\|^{2p} ){\rm sig}^{1-2p}\left( s_i \right) \notag \\
			&  \!+\! k {\rm sgn}\!\left(s_i\right) \!+\! \frac{2c}{p(1\!-\!\eta)T_m} \!\sum\limits_{j = 1}^N  \! \exp(\!(a_{ij}\!\left\| \!{{x_i} \!-\! {x_j}}\! \right\|^2\! )^p) \!\times \notag \\
			& a_{ij}^{1-p}{\rm sig}^{1-2p}\left( {x_i} - {x_j} \right)  \Big) + d_i(t), \\
			{s_i}(t) =&~ \!\nabla \!{f_i}\!\left( \!{{x_i}}\! \right) \!+\!\! \int_0^t \!\! \frac{2c}{p(1\!-\!\eta)T_m} \!\Big( \!\sum\limits_{j = 1}^N \! \exp(\!(a_{ij}\!\left\| {{x_i} - {x_j}} \right\|^2\! )^p) \! \times \notag \\
			& a_{ij}^{1-p}{\rm sig}^{1-2p}\left( {x_i} - {x_j} \right)  \Big) d\tau 
		\end{align}
	\end{subequations}
	where $ 0< \eta <1 $, $ 0 < p < {1}/{2}$, $T_m>0$, $ c \ge {\bar{\varPsi}}/{(4 \lambda _2 \left( \mathcal{L} \right))} $, and $ k \ge \mathcal{H}\mathcal{D} $ with $ \left\| {\nabla ^2}{f_i}\left( {{x_i}} \right) \right\|_1 \le \mathcal{H}$. The proof details can be done similarly to that of Theorem \ref{theorom_PT_Algorithm}. The same results can be found in \cite[Proposition 2]{Guo2023local-minimization-free}, where an additional sliding manifold is involved to deal with disturbance. In contrast, algorithm \eqref{PdT_Algorithm_Disturbance} can guarantee both zero-gradient-sum and disturbance-rejection by a common sliding mode.
\end{remark}

\begin{remark}[Full Distributed Implementation]
	Since the parameter selection of algorithm \eqref{Predefined_Algorithm} depends on the algebraic connectivity of the graph and maximum Lipschitz constant of local gradients, it is not fully distributed. However, the proposed framework can readily be extended to fully distributed approximate optimization with the help of the notion of time-base generators in \cite{Guo2020predefined-time,Liu2023multi-objective,Guo2022distributed_Auto}. A predefined-time approximate stability can be achieved by multiplying a well-defined dynamic system with a time-varying gain generated by a time-base generator. Consider the following algorithm
	%	With the help of the time-base generator used in \cite{Guo2020predefined-time,Liu2023multi-objective,Guo2022distributed_Auto,Li2022predefined-time}, we can achieve fully distributed predefined time approximation optimization, as shown in detail below:
	\begin{subequations}
		\label{PdT_Algorithm_Full}
		\begin{align}
			\dot{x}_i(t) =&~ - \big({ {\nabla ^2}{f_i}\left( {{x_i}} \right)}\!\big)^{-1} \!\Big(  \frac{1}{2p\eta T_m}\exp(\left\| s_i \right\|^{2p} ){\rm sig}^{1-2p}\left( s_i \right) \notag \\
			&+ \mathbb{T} \left(t,\left(1-\eta \right)T_m \right) \sum\limits_{j = 1}^N { {a_{ij}}\left( {{x_i} - {x_j}} \right) }  \Big), \\
			{s_i}(t) =&~ \nabla\! {f_i}\left( {{x_i}} \right) \!+\! \int_0^t \! \mathbb{T} \left(t,\left(1\!-\!\eta \right)T_m \right) \! \sum\limits_{j = 1}^N \!{ {a_{ij}}\left( {{x_i} \!-\! {x_j}} \right) } d\tau 
		\end{align}
	\end{subequations}
	where $ 0< \eta <1 $, $ 0 < p < {1}/{2}$, $T_m >0$, and $ \mathbb{T} \left(t,\left(1-\eta \right)T_m \right)  $ is a time-base generator as defined in \cite{Liu2023multi-objective}. \textcolor{black}{It is obvious that the zero-gradient-sum is achieved within $\eta T_m$ and the states of all agents to the small domain of the optimal solution within $T_m$ under algorithm \eqref{PdT_Algorithm_Full}. The proof can be done along similar lines as that of Theorem 1, which is omitted here to save space.}
\end{remark}

\begin{remark}
	It is worth noting that the communication topology considered in this paper is fixed, and the connection weights are assumed to be piecewise constant even in the dynamic switching topology. However, sensing or communication capability is range-limited, so it is impractical to simply make an assumption that network connectivity is preserved by default. One potential solution is to maintain network connectivity with artificial potential functions \cite{Ning2019TCYB,Hong2020TCNS,Wu2022TCYB}. In the future, we will investigate how to implement ZGS algorithms for network connectivity preservation under preset time constraints.
\end{remark}

	Finally, in the following, we summarize the proposed algorithm \eqref{Predefined_Algorithm} in Algorithm \ref{algorithm1}.
\begin{algorithm}[htb]
	\caption{Predefined-time Distributed Optimization Algorithm}
	\label{algorithm1}
	\begin{algorithmic}[1] % 1代表显示行号
		\REQUIRE ~~\\ %算法的输入：Input
		$x_i(0)$,  $ \nabla f_i $,  $ \nabla ^2 f_i$, $T_m$, $\bar{\varPsi}$, $\lambda _2 \left( \mathcal{L} \right)$;
		\ENSURE ~~\\ %算法的输出：Output
		$x^*$;
		\STATE \textbf{Initialization:}  \\
		$t=0$;\\
		Select parameters $\eta,p,c$ according to Theorem 1;\\
		\WHILE{$t\le T_m$}
		\FOR {each agent $i \in \cal{V}$ }
		\STATE Update $x(t)$ by Eq. \eqref{Predefined_Algorithm}; \\
		\ENDFOR
		\STATE Update $t$; \\
		\ENDWHILE
		\RETURN $x^* =x_i(T_m)$;
	\end{algorithmic}
\end{algorithm}

\subsection{Predefined-time Distributed Time-varying Optimization Algorithms}
Most existing distributed time-varying optimization methods are distributed implementations of prediction-correction algorithms utilizing average tracking technology (see \cite{Salar2017distributed} and its extensions). The time-varying optimization with predefined-time convergence under the ZGS framework is developed in this section.
\begin{assumption}
	\label{assumption_cost_time-varying}
	Each cost function $f_i$ has the identical Hessian and the rate of change of local gradients is bounded, i.e., $ \nabla^2 f_i\left(  x_i ,t \right) = \nabla^2 f_j\left(  x_j ,t \right) $, $ \| \frac{\partial}{\partial t}\nabla f_i\left(  x_i ,t \right) \|_\infty \le \omega$ for $i,j\in \mathcal{V}$ and $t\ge 0$.
\end{assumption}
\begin{remark}
	Assumption \ref{assumption_cost_time-varying} is trivial and often considered in time-varying optimization problems. It is easy to find functions that satisfy the requirements of Assumption \ref{assumption_cost_time-varying}, such as quadratic functions $f(x,t) = (\gamma x +\theta (t))^2$ commonly used for energy optimization and engineering application.
\end{remark}

For each agent $i$, the distributed algorithm is designed as follows
\begin{subequations}
	\label{Prescribed_Algorithm_time-varying}
	\begin{align}
		\label{PT_Algorithm_dynamic_time-varying}
		\dot{x}_i(t) =& - \big({ {\nabla ^2}{f_i}\left( {{x_i}} \right)}\big)^{-1} \Big(  \frac{1}{2p\eta T_m}\exp(\left\| s_i \right\|^{2p} ) {\rm sig}^{1-2p}\left( s_i \right) +  \notag \\
		& \frac{2c}{p(1\!-\!\eta)T_m} \!\sum\limits_{j = 1}^N {\! \exp\!\big(\!(a_{ij}\left\| \!{{x_i} \!-\! {x_j}} \!\right\|^2\! )^p\!\big)\! a_{ij}^{1\!-\!p}\!\textcolor{blue}{{\rm SIG}}^{1\!-\!2p}\!\left(\! {x_i} \!-\! {x_j}\! \right) }  \notag \\
		&+ \mu \sum\limits_{j = 1}^N {{a_{ij}}\! {\rm \textcolor{blue}{SGN}} \left( {{x_i} - {x_j}} \right)} + \frac{\partial }{\partial t} \nabla f_i\left(  x_i ,t \right) \Big), \\
		\label{PT_Algorithm_sliding mode_time-varying}
		{s_i}(t) =& \nabla \! {f_i}\left( {{x_i}} \right) \!+\! \int_0^t \! \Big(\!\mu \!\sum\limits_{j = 1}^N\! {{a_{ij}} {\rm \textcolor{blue}{SGN}} \left( {{x_i} \!-\! {x_j}} \right)} \!+\! \frac{2c}{p(1\!-\!\eta)T_m}\! \times \notag \\
		&\sum\limits_{j = 1}^N { \!\exp\!\big(\!(a_{ij}\left\| {{x_i}\! -\! {x_j}} \right\|^2 )^p\!\big) a_{ij}^{1-p}\textcolor{blue}{\rm SIG}^{1-2p}\left( {x_i} \!-\! {x_j} \right) } \Big) d\tau 
	\end{align}
\end{subequations}
where $ \eta, p, T_m, c, \mu \in \mathbb{R}$ are positive constants.

\begin{theorem}
	\label{theorom_PT_Algorithm_time-varying}
	Suppose that Assumptions \ref{assumption_G} - \ref{assumption_cost_time-varying} hold. Let $ \underline{\psi}= \min_{ i \in {\mathcal V} } \{\psi _i\} $, $ \bar{\varPsi} = \max_{ i \in {\mathcal V} } \{\varPsi _i\} $, and $ \underline{a} = \min_{i,j \in \mathcal{V}}\{a_{ij} :a_{ij}\neq 0 \}$. If $ 0< \eta <1 $, $ 0 < p < {1}/{2}$, $ T_m >0 $, $ c \ge {\bar{\varPsi}}/\left({4 \lambda _2 \left( \mathcal{L} \right)}\right) $, and $ \mu \ge {2N \omega \bar{\varPsi} }/{(\underline{a} \underline{\psi})}$, then, distributed algorithm \eqref{Prescribed_Algorithm_time-varying} drives the states $x_i(t)$ of all agents to a optimal solution $x^*(t)$ of \eqref{optimization_problem_Time-varying} within a settling time that is bounded by $ T_m $.
\end{theorem}

\begin{IEEEproof}For $ 0<t \le \eta T_m$, differentiating \eqref{PT_Algorithm_sliding mode_time-varying} 
yields
\begin{align}
	\label{PT_Algorithm_sliding mode_dynamic_time-varying}
	&\dot s _i (t) = \nabla^2 f_i \left( x_i,t \right) \dot x_i \!+\! \frac{\partial }{\partial t} \nabla\! f_i\left(  x_i ,t \right) \!+\! \mu \!\sum\limits_{j = 1}^N {{a_{ij}} {\rm \textcolor{blue}{SGN}} \left( {{x_i} \!-\! {x_j}} \right)} \notag \\
	&\!+\! \frac{2c}{p(1\!-\!\eta)T_m} \!\sum\limits_{j = 1}^N { \!\exp\!\big(\!(a_{ij}\left\| {{x_i} \!-\! {x_j}} \right\|^2 )^p\!\big) \!a_{ij}^{1\!-\!p}\!\textcolor{blue}{\rm SIG}^{1\!-\!2p}\left({x_i} \!-\! {x_j} \right) } \notag \\
	&= - \frac{1}{2p\eta T_m}\exp(\left\| s_i \right\|^{2p} ){\rm sig}^{1-2p}\left( s_i \right)
\end{align}
Similar to Theorem \ref{theorom_PT_Algorithm}, it holds that $ \lim_{t \to \eta T_m} s_i \left( t \right) = {\bf 0} $.

For $ t > \eta T_m$, the dynamic satisfies 
\begin{align}
	\label{PT_Algorithm_closed_loop_dynamic_time-varying}
	&\dot{x}_i (t) = - \!\big({ {\nabla ^2}\!{f_i}\left( {{x_i}},t \right)}\!\big)^{-1} \!\Big(\! \mu\! \sum\limits_{j = 1}^N {\!{a_{ij}} {\rm \textcolor{blue}{SGN}} \left( {{x_i} \!-\! {x_j}} \right)} \!+\! \frac{\partial }{\partial t} \nabla \!f_i\left( \! x_i ,t \right) \notag \\
	&+ \!\frac{2c}{p(1\!-\!\eta)T_m} \!\sum\limits_{j = 1}^N { \!\exp\!\big(\!(a_{ij}\!\left\| {{x_i} \!-\! {x_j}} \right\|^2 )^p\!\big) a_{ij}^{1\!-\!p}\textcolor{blue}{\rm SIG}^{1\!-\!2p}\left( {x_i} \!-\! {x_j} \right) } \!\Big) 
\end{align}
It follows from \eqref{PT_Algorithm_closed_loop_dynamic_time-varying} that ${{d \sum\nolimits_{i = 1}^N {\nabla {f_i}\left( {{x_i},t} \right)}}}/{{dt}} = \sum\nolimits_{i = 1}^N ({\nabla^2 {f_i}\left( {{x_i},t} \right)} \dot{x}_i \!+\! \frac{\partial }{\partial t} \nabla f_i\left(  x_i ,t \right) )  = {\bf 0} $. Due to $ \sum\nolimits_{i = 1}^N {\nabla {f_i}\left( {{x_i},t} \right)}  = {\bf 0} $ at $ t = \eta T_m $, one has $ \sum\nolimits_{i = 1}^N {\nabla {f_i}\left( {{x_i},t} \right)}  = {\bf 0} \forall t \ge \eta T_m $. 

Denote $ H'(t) = { {\nabla ^2}{f_i}\left( {{x_i}},t \right)}^{-1}$ and $ {\bf \tilde x} = \left[\tilde x_1^T,\ldots, \tilde x_N^T \right]^T \in \mathbb{R}^{nN} $, where ${\tilde x_i} = x_i - {\bar x} $ with ${\bar x} = \frac{1}{N}\sum_{k=1}^{N}x_k $. Consider the following Lyapunov function for system (\ref{PT_Algorithm_closed_loop_dynamic_time-varying})
\begin{align}
	{W} \left(x_i(t)\right) = \frac{1}{2} \left\| {\bf \tilde x} \right\|^2 = \frac{1}{2}\sum_{i=1}^{N} \left\| {\tilde x_i} \right\|^2 
\end{align}
Differentiating $ W  $ along (\ref{PT_Algorithm_closed_loop_dynamic_time-varying}) yields
\begin{align}
	\dot{W} =& -\sum_{i=1}^{N} {\tilde x_i} ^T \bigg( H'(t) \Big( \mu \sum\limits_{j = 1}^N {{a_{ij}} {\rm \textcolor{blue}{SGN}} \left( {{x_i} - {x_j}} \right)} \!+\! \frac{\partial }{\partial t} \nabla f_i\left(  x_i ,t \right) \notag \\
	&+\frac{2c}{p(1-\eta)T_m} \sum\limits_{j = 1}^N  \exp\big((a_{ij}\left\| {{x_i} - {x_j}} \right\|^2 )^p\big) \times \notag \\
	&a_{ij}^{1-p}\textcolor{blue}{\rm SIG}^{1-2p}\left( {x_i} - {x_j} \right)  \Big) \bigg) \notag \\
	=& - \!\mu \!\sum_{i=1}^{N} {\tilde x_i} ^T  H'(t) \big(\!\sum\limits_{j = 1}^N {{a_{ij}} {\rm \textcolor{blue}{SGN}} \left( {{\tilde x_i} - {\tilde x_j}} \right)}\! +\! \frac{\partial }{\partial t} \nabla f_i\left(  x_i ,t \right)\! \big) \notag \\
	&- \!\frac{2c}{p(1-\eta)T_m} \!\sum_{i=1}^{N} {\tilde x_i} ^T  H'(t)  \!\sum\limits_{j = 1}^N  \exp\big(\!(a_{ij}\left\| {{\tilde x_i} \!- \!{\tilde x_j}} \right\|^2 )^p\!\big) \times \notag \\
	&a_{ij}^{1-p}\textcolor{blue}{\rm SIG}^{1-2p}\left( {\tilde x_i} - {\tilde x_j} \right)  \notag \\
	\le&-\frac{\mu }{2\bar{\varPsi}} \sum_{i=1}^{N} \sum\limits_{j = 1}^N {{a_{ij}} \textcolor{blue}{\left\| {{\tilde x_i} - {\tilde x_j}} \right\|} }  + \frac{\omega}{{\underline{\psi}}}\sum_{i=1}^{N} { \left\| {\tilde x_i}  \right\|} \notag \\
	& - \frac{c}{{\bar{\varPsi}} p(1-\eta)T_m} \sum\limits_{i = 1}^N \sum\limits_{j = 1}^N  \Big( \exp\big((a_{ij}\left\| {{\tilde x_i} - {\tilde x_j}} \right\|^2 )^p\big)  \times \notag \\
	&  a_{ij}^{1-p} \textcolor{blue}{\left\| \tilde x_i - \tilde x_j \right\|^{2-2p}} \Big) 
\end{align}
Notice that \textcolor{blue}{
	$ { \left\| {\tilde x_i}  \right\|} \le \frac{1}{N} \sum_{j=1}^{N} { \left\| {x_i}-{x_j}  \right\|} \le \sum_{i=1}^{N} \sum_{j=1}^{N} { \left\| {x_i}-{x_j}  \right\|} = \sum_{i=1}^{N} \sum_{j=1}^{N} { \left\| {\tilde x_i} - {\tilde x_j}  \right\|} $.}
%Notice that 
%$ { \left\| {\tilde x_i}  \right\|} \le \frac{1}{N} \sum_{j=1}^{N} { \left\| {x_i}-{x_j}  \right\|} \le \sum_{i=1}^{N} \sum_{j=1}^{N} { \left\| {x_i}-{x_j}  \right\|} = \sum_{i=1}^{N} \sum_{j=1}^{N} { \left\| {\tilde x_i} - {\tilde x_j}  \right\|} \le \sum_{i=1}^{N} \sum_{j=1}^{N} { \left\| {\tilde x_i} - {\tilde x_j}  \right\|_1} $.
It can be concluded that
\begin{align}
	\dot{W}_{C} \le& -\frac{\mu }{2\bar{\varPsi}} \sum_{i=1}^{N} \sum\limits_{j = 1}^N {{a_{ij}} \textcolor{blue}{\left\| {{\tilde x_i} - {\tilde x_j}} \right\|}} \!+\! \frac{N \omega}{{\underline{\psi}}}\sum_{i=1}^{N} \sum_{j=1}^{N} \textcolor{blue}{ \left\| {\tilde x_i} - {\tilde x_j}  \right\|} \notag \\
	& - \frac{c}{{\bar{\varPsi}} p(1-\eta)T_m} \sum\limits_{i = 1}^N \sum\limits_{j = 1}^N  \Big( \exp\big((a_{ij}\left\| {{\tilde x_i} - {\tilde x_j}} \right\|^2 )^p\big)  \times \notag \\
	&  { \big( {a_{ij}} \left\| {{\tilde x_i} - {\tilde x_j}} \right\|^2 \big)^{1-p} } \Big) 
\end{align}
Due to $ \mu \ge {2N\omega \bar{\varPsi} }/(\underline{a}{\underline{\psi}})$ and $ c \ge {\bar{\varPsi}}/\left({4 \lambda _2 \left( \mathcal{L} \right)}\right) $, it is clear that 
\begin{align}
	\dot{W} \le& -\! \frac{1}{p(1\!-\!\eta)T_m} \!\big(\! \frac{ 4{\lambda _2 \left( \mathcal{L} \right)} }{N(N\!-\!1)} \!\big)^{-p} \!\exp\!\big(\! ( \frac{4{\lambda _2 \left( \mathcal{L} \right)} W }{N(N\!-\!1)} )^p\big)  { W^{1-p} }
\end{align} 
By Lemma \ref{lemma_PdT-time} and optimality conditions, the conclusion can be obtained.
\end{IEEEproof}

\begin{remark}
		Note that $T_m$ is an upper bound on the settling time of the algorithms \eqref{Predefined_Algorithm} and \eqref{Prescribed_Algorithm_time-varying}. Given the desired convergence time $T_m$, the different performance of the algorithms can be obtained by adjusting other parameters $ \eta, p, c$ and $\mu$. However, to be used for physical systems, the gains cannot be too large due to actuator saturation. Therefore, the choice of controller parameters is a tradeoff between optimization performance and practical limitations. Besides,  the discontinuous function $\mathop{\rm \textcolor{blue}{SGN}}\left(\cdot\right)$ in \eqref{Prescribed_Algorithm_time-varying} may lead to the chattering of control input. To reduce chattering, one can replace $\mathop{\rm \textcolor{blue}{SGN}}\left(z\right)$ with a continuous approximation $\hat{\rm \textcolor{blue}{SGN}}(z)=\frac{z}{\|z\|+\epsilon(t)}$, where $ \dot \epsilon(t) = - \frac{1}{2p\eta T_m}\exp(\left\| \epsilon \right\|^{2p} ){\rm sig}^{1-2p}\left( \epsilon \right) $, which will reduce the chattering effect and make the controller easier to implement in real applications\cite{Zhao2022TAC-Design}.
\end{remark}

	The proposed time-varying algorithm \eqref{Prescribed_Algorithm_time-varying} is summarized in Algorithm \ref{algorithm2}.
\begin{algorithm}[htb]
	\caption{Predefined-time Distributed Time-varying Optimization Algorithm}
	\label{algorithm2}
	\begin{algorithmic}[1] % 1代表显示行号
		\REQUIRE ~~\\ %算法的输入：Input
		$x_i(0)$,  $ \nabla f_i $,  $ \nabla ^2 f_i$, $\frac{\partial}{\partial t}\nabla f_i$, $T_m$, $\omega$, $\underline{\psi}$, $\bar{\varPsi}$, $\lambda _2 \left( \mathcal{L} \right)$, $\underline{a}$;
		\ENSURE ~~\\ %算法的输出：Output
		$x^*(t)$;
		\STATE \textbf{Initialization:}  \\
		$t=0$;\\
		Select parameters $\eta,p,c,\mu$ according to Theorem 2;\\
		\WHILE{$t\ge 0$}
		\FOR {each agent $i \in \cal{V}$ }
		\STATE Update $x(t)$ by Eq. \eqref{Prescribed_Algorithm_time-varying}; \\
		\ENDFOR
		\STATE Update $t$; \\
		\IF{$t\ge T_m$}
		\RETURN $x^*(t) =x_i(t)$;
		\ENDIF
		\ENDWHILE
	\end{algorithmic}
\end{algorithm}

\begin{figure}[htb]
	\color{blue}
	\begin{center}
		\includegraphics[width=3.5cm]{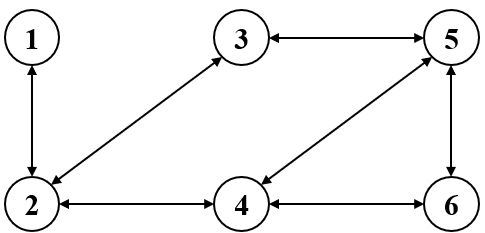} 
		\caption{Communication topology for multiagent systems (all weights are 1).}   
		\label{fig_communication}                              
	\end{center}                               
\end{figure}

\section{Simulation Verification}
In this section, the two simulation examples are provided to verify distributed time-invariant algorithms \eqref{Predefined_Algorithm} and time-varying algorithms \eqref{Prescribed_Algorithm_time-varying}, respectively. All simulations are performed using MATLAB/Simulink with an Euler solver and fundamental step-size $10^{-3}$.

\subsection{A Numerical Case}
\begin{figure}[htb]
	\begin{center}
		\subfloat[]{ \includegraphics[width=5.5cm]{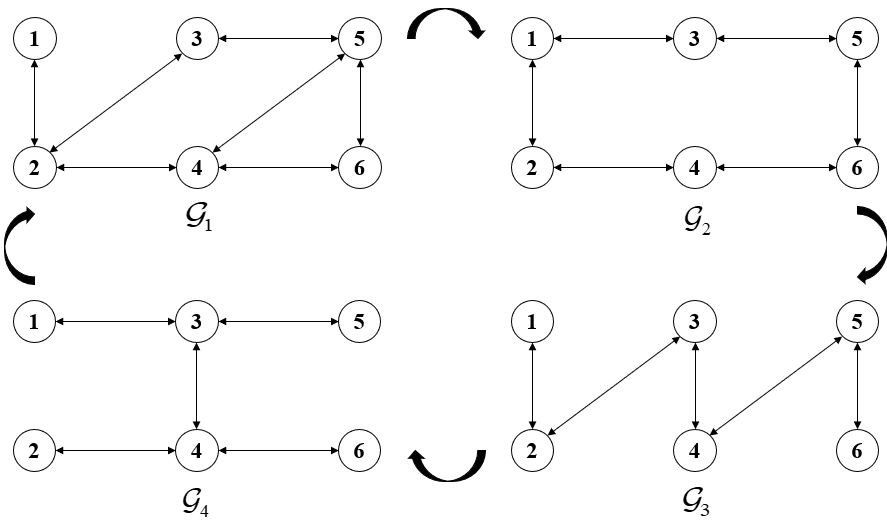}  \label{fig_communication_a_TV} }
		\hfil
		\subfloat[]{ \includegraphics[width=5cm]{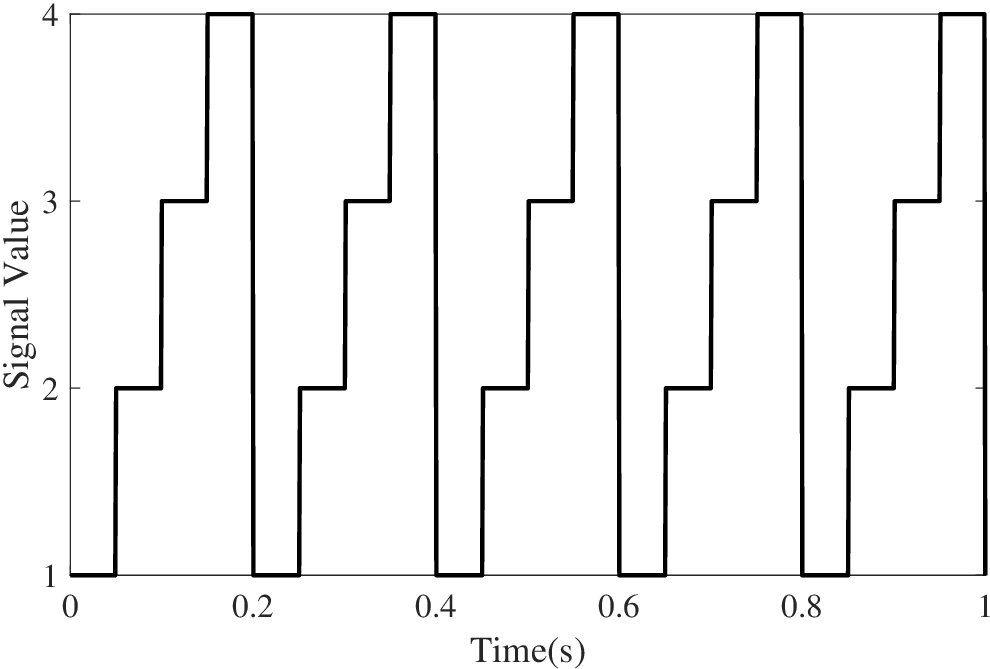} \label{fig_communication_a_signal}  }	
		\caption{Switching communication topologies for multiagent systems. (a) Set of graphs. (b) Switching signal.}                                 
	\end{center}                              
\end{figure}

\begin{figure}[htb]
	\centering
	\subfloat[]{ \includegraphics[width=4.25cm]{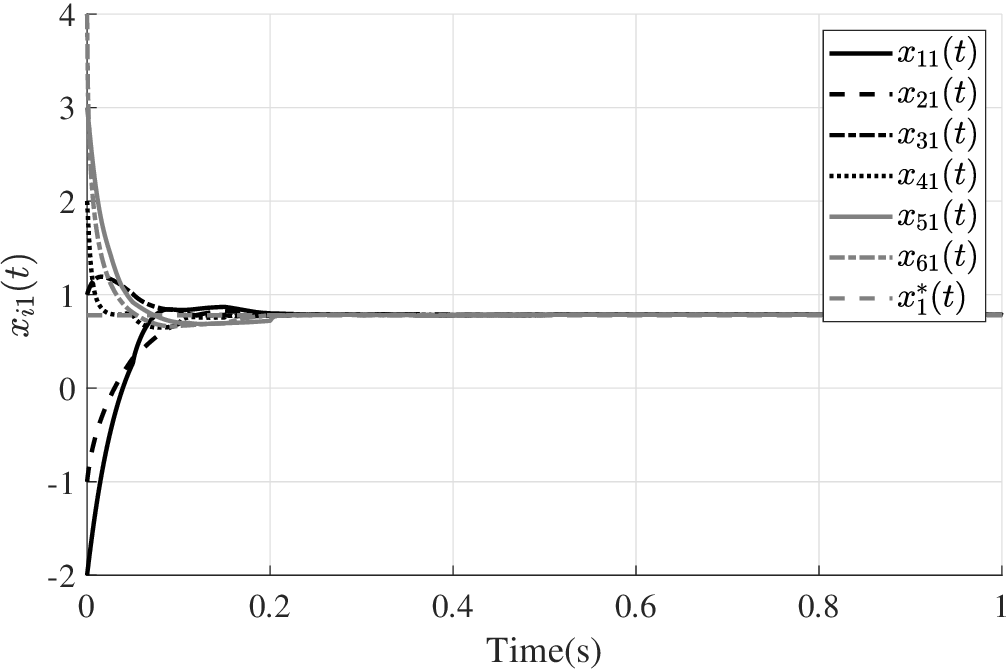} }
%	\hfil 
	\subfloat[]{ \includegraphics[width=4.2cm]{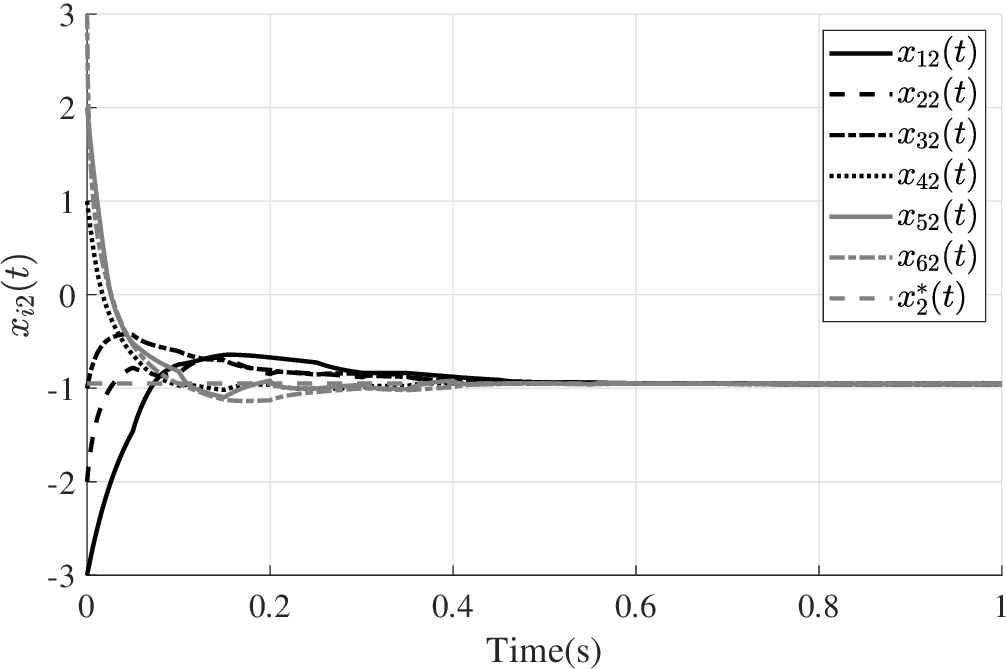} }	
	\caption{Evolution of $ x_i\left(t\right) $ under algorithm (\ref{Predefined_Algorithm} over switching graph). (a) $ x_{i1}\left(t\right) $. (b) $ x_{i2}\left(t\right) $ }
	\label{fig-PdT_LMFZGS_TV}
\end{figure}

Firstly, let us verify the distributed time-invariant optimization algorithm \eqref{Predefined_Algorithm}.
Consider a multi-agent system with 6 nodes to solve the global optimization problem $\min \sum_{i=1}^{6}f_i$ cooperatively. The local cost functions $f_i$ assigned to agents  are listed as: 
\begin{align*}
	f_1(x) &= (x_1 - 0.5)^2+2(x_2 +1.3)^2 -0.5x_1x_2 ,\\
	f_2(x) &= 2(x_1+0.7)^2 +1.5(x_2+1.7)^2 +0.3x_1x_2 \\
	&+0.3\sin (0.3x_1+1.8)+0.73\cos(0.5x_2+1) ,\\
	f_3(x) &= 2(x_1-1.5)^2 +2(x_2-0.3)^2 + \ln(2+0.1x_1^2) \\
	&+ \ln (4+0.6x_2^2), \\
	f_4(x) &= 0.5(x_1-1.5)^2 +1.5(x_2+1.6)^2 +0.5x_1x_2 \\
	&+ x_1/(\sqrt{2+0.4x_1^2}) + 0.6x_2/(\sqrt{1+x_2^2}), \\
	f_5(x) &= (x_1-2)^2 +(x_2-0.9)^2 +0.7x_1x_2 \\
	&+ 0.3\exp(-0.4x_1^2) + 0.7\exp(-0.5x_2^2) , \\
	f_6(x) &= 1.5(x_1-0.8)^2+2(x_2+1.5)^2.
\end{align*}
  \begin{figure}[!t]
 	\centering
 	\subfloat[]{ \includegraphics[width=4.25cm]{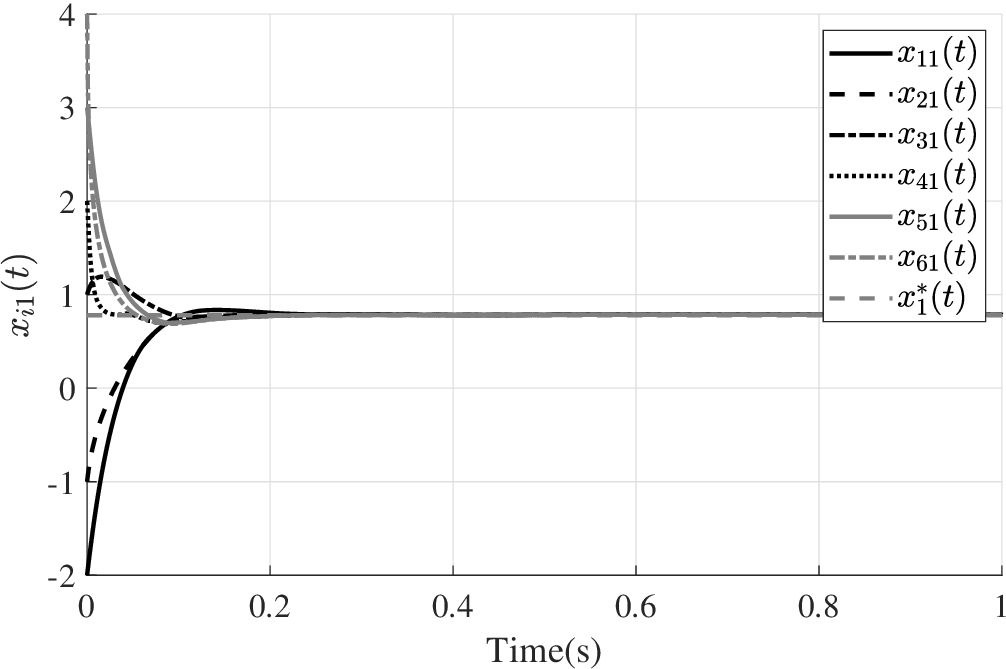} } %width=\columnwidth
%  	\hfil 
 	\subfloat[]{ \includegraphics[width=4.2cm]{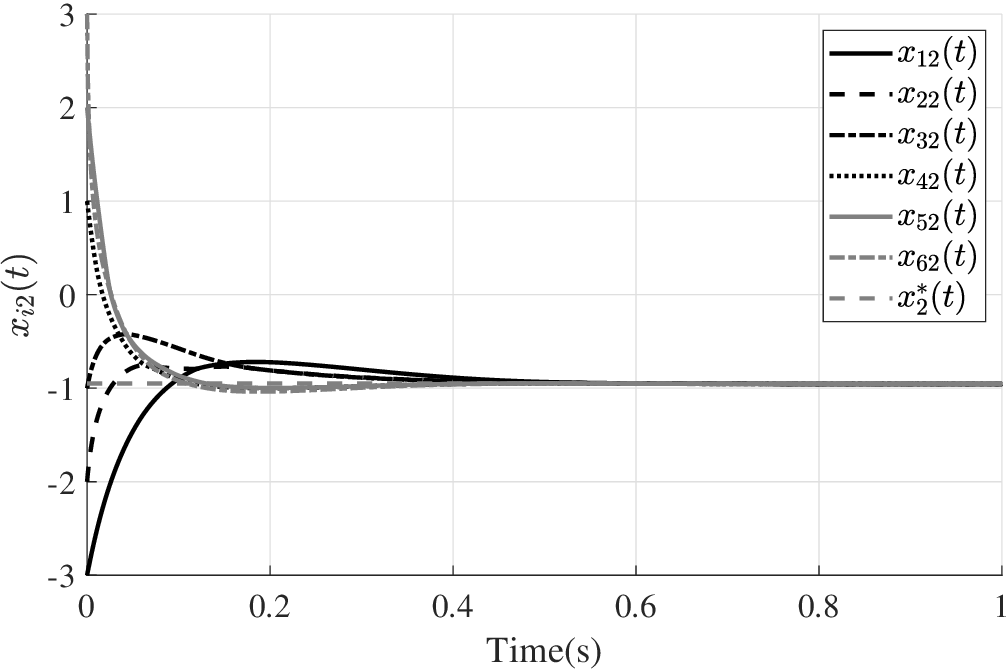} }	
 	\caption{Evolution of $ x_i\left(t\right) $ under algorithm (\ref{Predefined_Algorithm}). (a) $ x_{i1}\left(t\right) $. (b) $ x_{i2}\left(t\right) $ }
 	\label{fig-PdT_LMFZGS}
 \end{figure}
 \begin{figure}[htb]
	\color{blue}
	\centering
	\includegraphics[width=6.3cm]{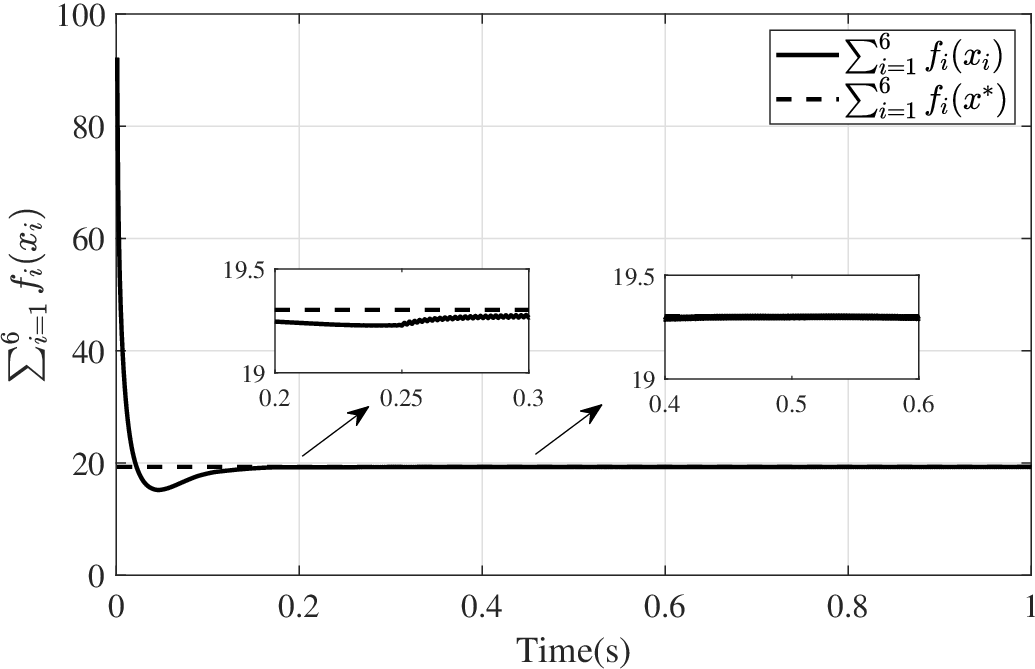}
	\caption{Evolution of the global objective function under algorithm (\ref{Predefined_Algorithm}).}
	\label{A3}
\end{figure}
 \begin{figure}[htb]
	\color{blue}
	\centering
	\subfloat{ \includegraphics[width=4.2cm]{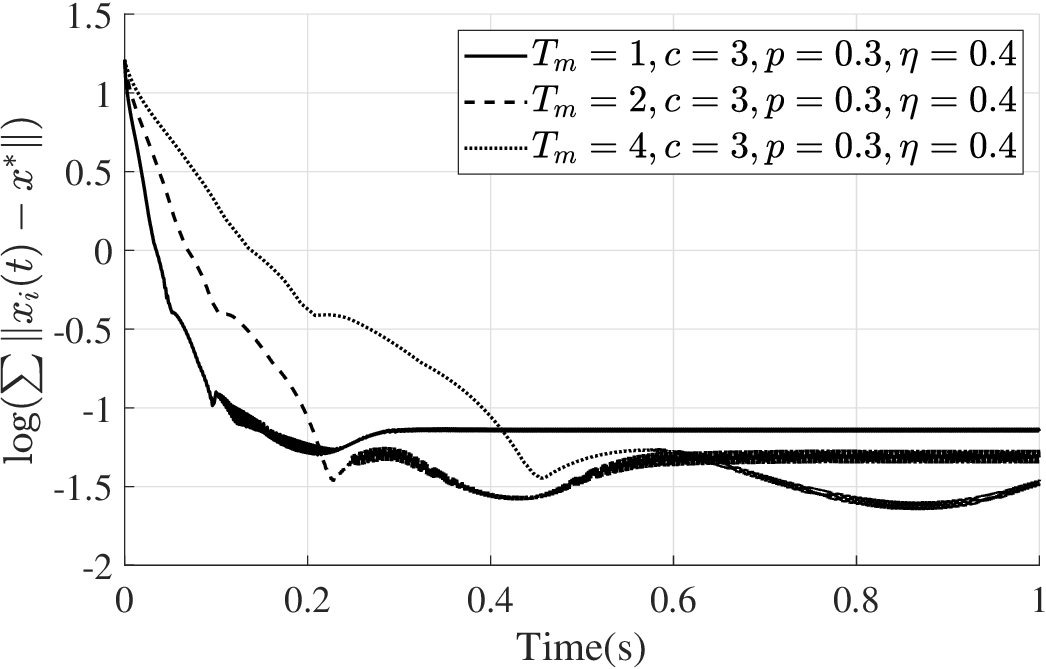} } %width=\columnwidth
	\subfloat{ \includegraphics[width=4.2cm]{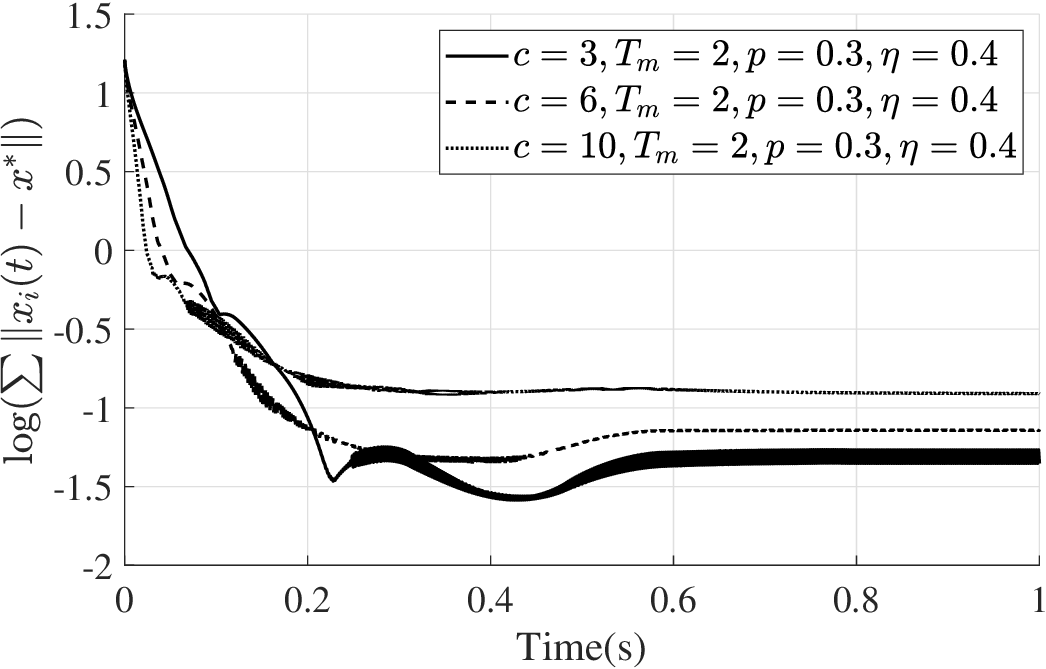} }	
	\hfil 
	\subfloat{ \includegraphics[width=4.2cm]{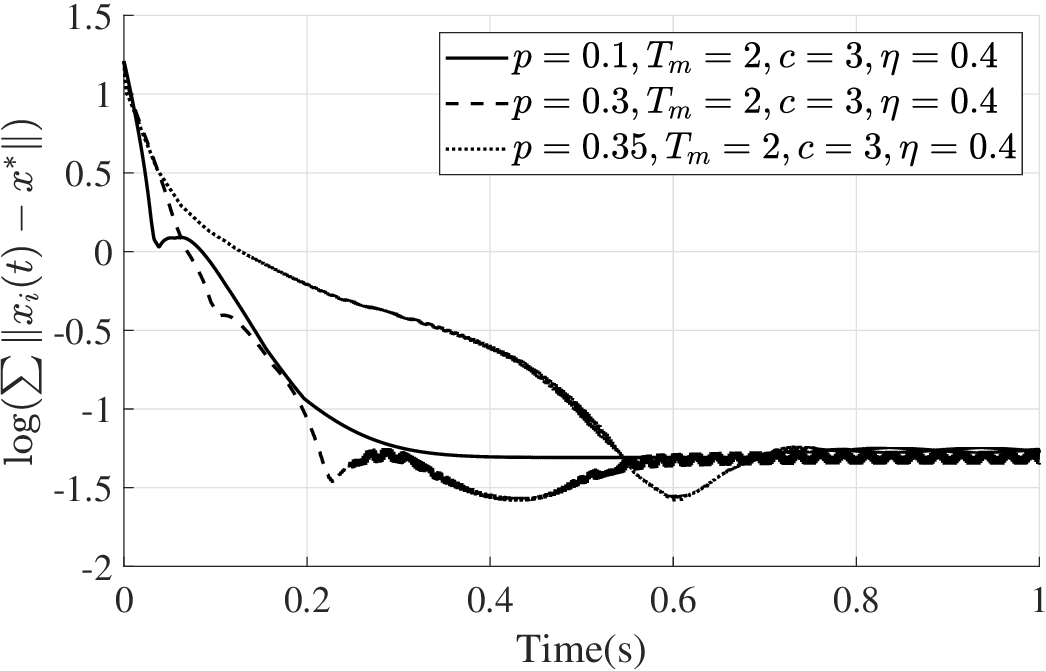} } %width=\columnwidth
	\subfloat{ \includegraphics[width=4.2cm]{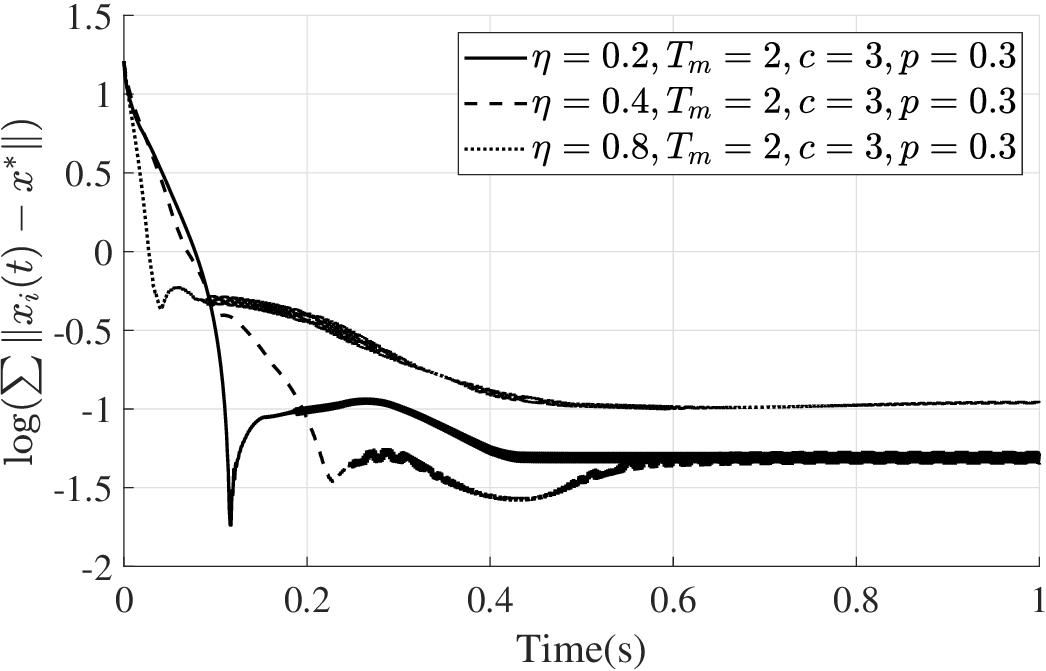} }	
	\caption{Performance evaluation of the algorithm (\ref{Predefined_Algorithm}) with different parameters.}
	\label{fig-PdT_LMFZGS_different_parameters}
\end{figure}
It can be verified that the optimal solution of global objective is $ x ^* =\left[x_1^*,x_2^* \right]^T= \left[0.7858,-0.9551 \right]^T $. The communication topology among agents is depicted as an undirected connected graph shown in Fig. \ref{fig_communication} and agents’ initial values are set as $ \boldsymbol{x}_{1}\left(0\right) =  \left[x_{i1}\left(0\right)\right]^T = \left[-2,-1,1,2,3,4\right]^T $ and $ \boldsymbol{x}_{2}\left(0\right) =  \left[x_{i2}\left(0\right)\right]^T = \left[-3,-2,-1,1,2,3\right]^T $. The algorithm parameters are selected as $ p=0.3 $, $ \eta = 0.4 $, $c= 3$ and $ T_m = 2$. Fig. \ref{fig-PdT_LMFZGS} depicts the evolution of the system states, from which it can be seen that the global optimizer performs accurately and
 quickly. The trajectory of the global objective function shown in Fig. \ref{A3} also illustrates the validity of algorithm \eqref{Predefined_Algorithm}. Fig. \ref{fig-PdT_LMFZGS_different_parameters} depicts the convergence performance of the algorithm with different parameters. As discussed in Remark 10, $T_m$ is an upper bound on the desired settling time, and a variety of convergence rates can be obtained by adjusting the parameters.
 Simulation of the algorithms under dynamic switching topologies is also done. Fig. \ref{fig_communication_a_TV} and \ref{fig_communication_a_signal} show the  the sequence and the switching signal of the topologies, respectively. The convergence results are depicted in Fig. \ref{fig-PdT_LMFZGS_TV}, which verifies the effectiveness of the algorithm.

\begin{figure*}[htb]
	\centering
	\subfloat[]{ \includegraphics[width=4.5cm]{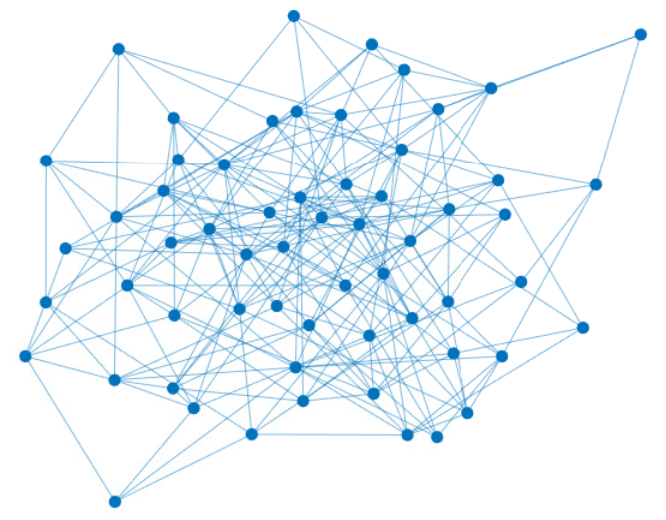} \label{fig-PdT_LMFZGS_60agent_a}}
	\subfloat[]{ \includegraphics[width=5.5cm]{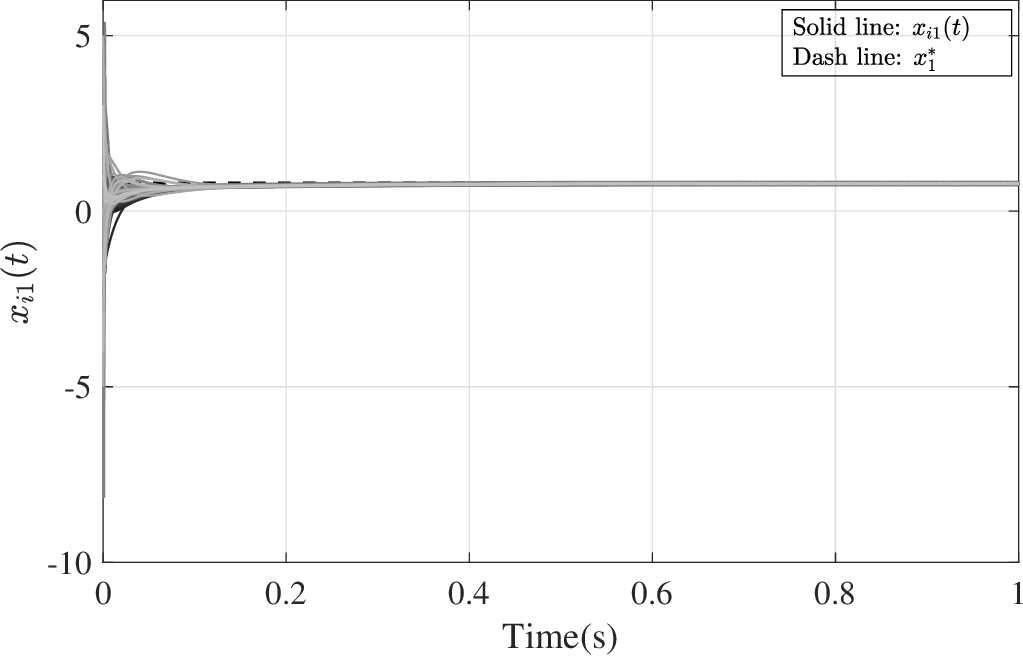} \label{fig-PdT_LMFZGS_60agent_b}} %width=\columnwidth
%	\hfil 
	\subfloat[]{ \includegraphics[width=5.4cm]{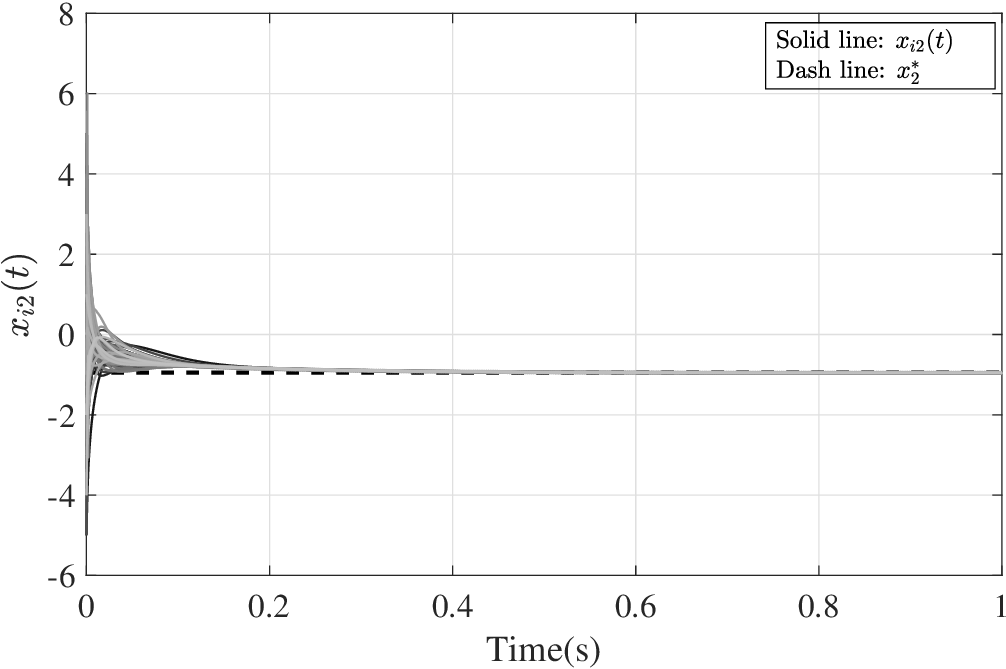} \label{fig-PdT_LMFZGS_60agent_c}}	
	\caption{Experiment of algorithm (\ref{Predefined_Algorithm}) over $60$ agents. (a) Communication topology. (b) Evolution of $ x_{i1}\left(t\right) $. (c) Evolution of $ x_{i2}\left(t\right) $. }
	\label{fig-PdT_LMFZGS_60agent}
\end{figure*}

To illustrate the advantages and scalability of the proposed distributed method, the algorithm is executed in a computing network consisting of $60$ agents. The communication topology is represented by a randomly generated connected graph, as shown in Fig. \ref{fig-PdT_LMFZGS_60agent_a}. The cost functions of agents are set to $f_{i*10+k} (x)= f_{i+1} (x), i = \{0,1,\dots,5\}, k = \{1,2,\dots,10\} $, where $f_1 $ - $  f_6$ are as previously determined. The algorithm parameters remain unchanged, and the initial values are randomly selected at $[-5,5]$. The convergence results are depicted in Fig. \ref{fig-PdT_LMFZGS_60agent_b} - \ref{fig-PdT_LMFZGS_60agent_c}, showing that the computational burden of the algorithm is not affected by the agent scale.

\begin{figure}[htb]
	\color{blue}
	\begin{center}
		\includegraphics[width=3.5cm]{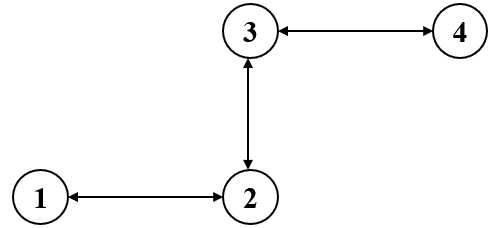} 
		\caption{Communication topology for multi-robots.} 
		\label{fig_communication_b}                                
	\end{center}                               
\end{figure}

\begin{figure}[htb]
	\centerline{\includegraphics[width=6cm]{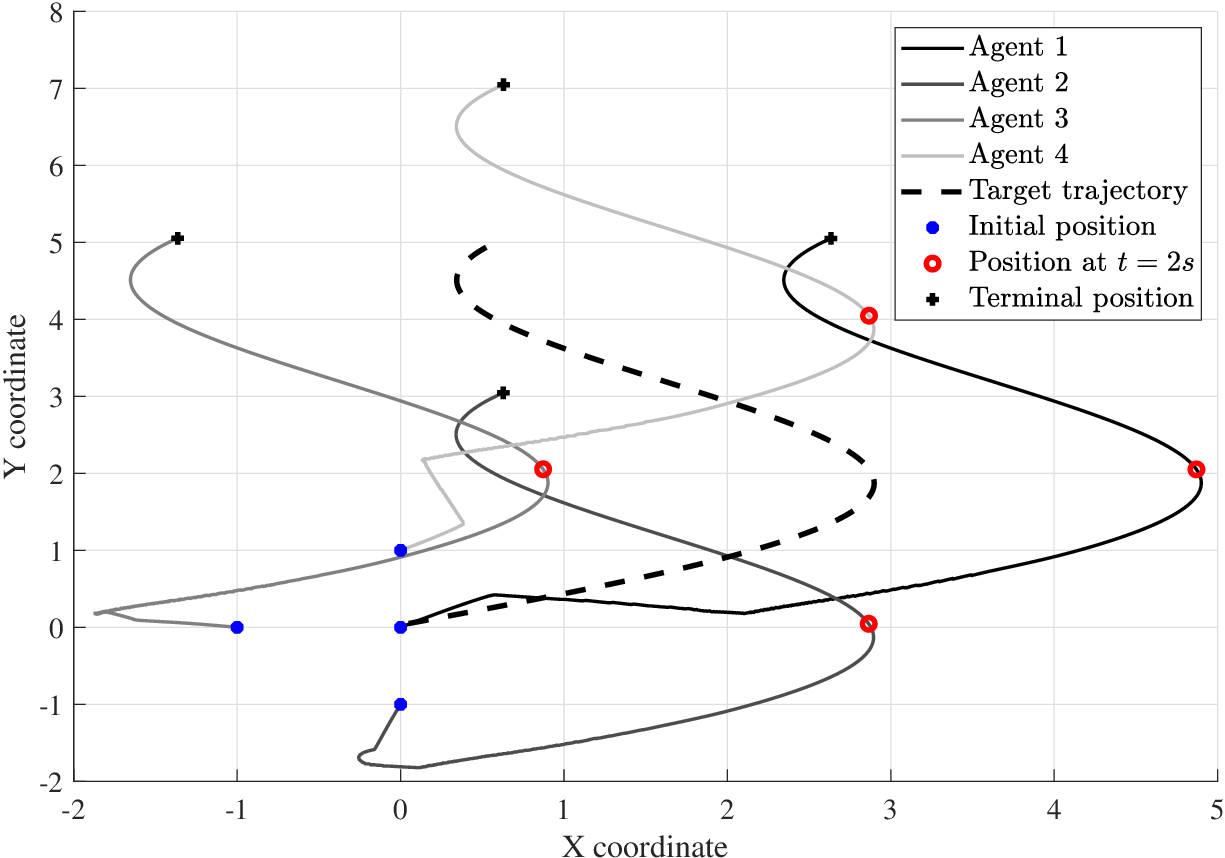}}
	\caption{Trajectories for robot $i$, $i \in \mathcal{V}$ and the target source.}
	\label{B1}
\end{figure}
\begin{figure}[htb]
	\centering
	\subfloat[]{ \includegraphics[width=4.25cm]{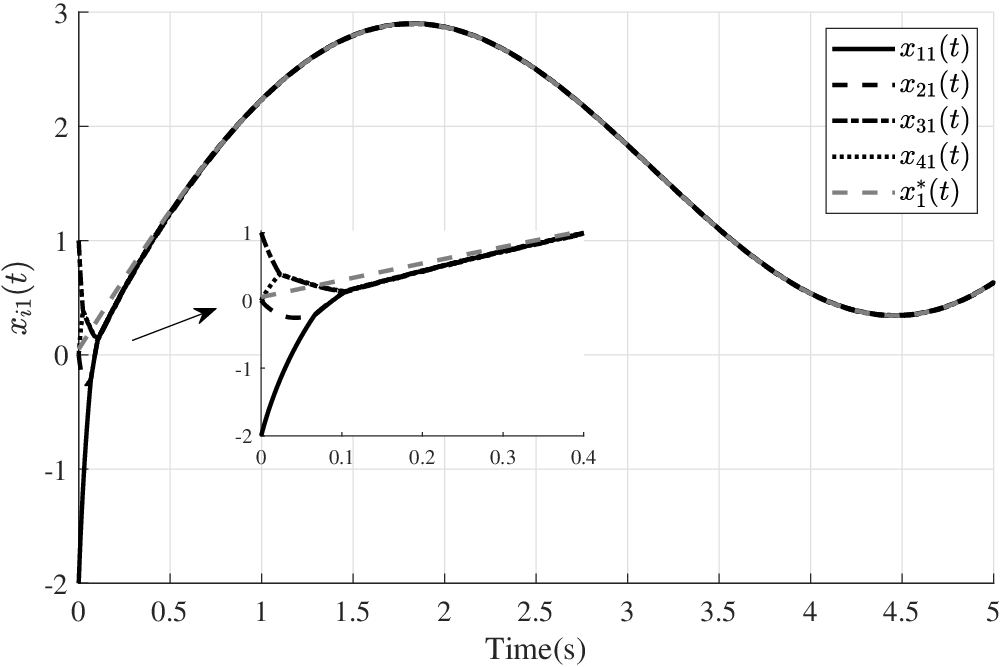} }
	%	\hfil 
	\subfloat[]{ \includegraphics[width=4.2cm]{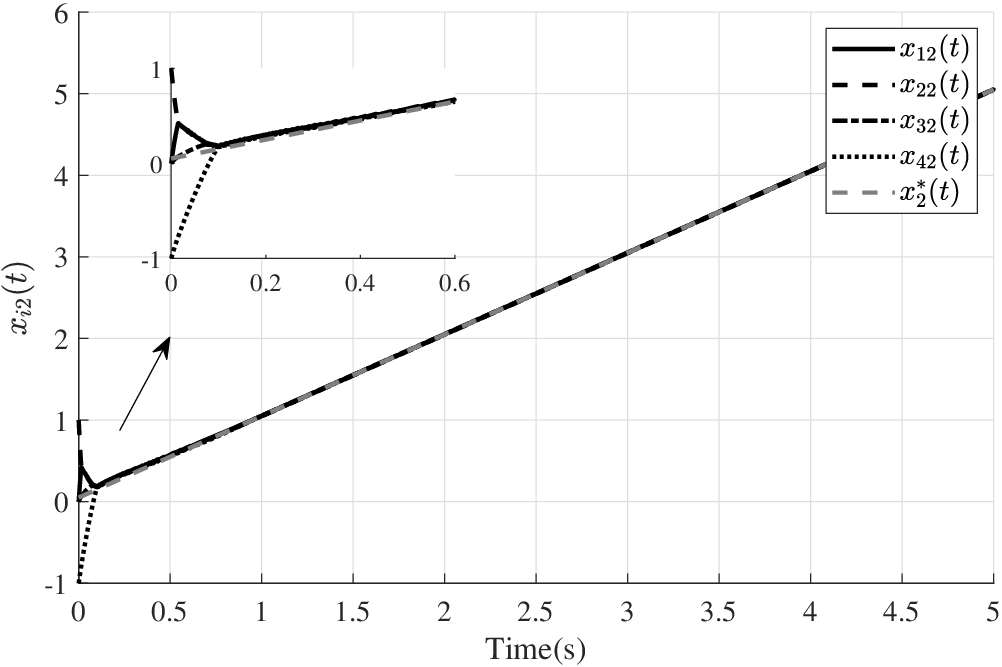} }	
	\caption{Evolution of $ x_i\left(t\right) $ under algorithm \eqref{Prescribed_Algorithm_time-varying}. (a) $ x_{i1}\left(t\right) $. (b) $ x_{i2}\left(t\right) $ }
	\label{B2}
\end{figure}
\begin{figure}[htb]
	\centerline{\includegraphics[width=6cm]{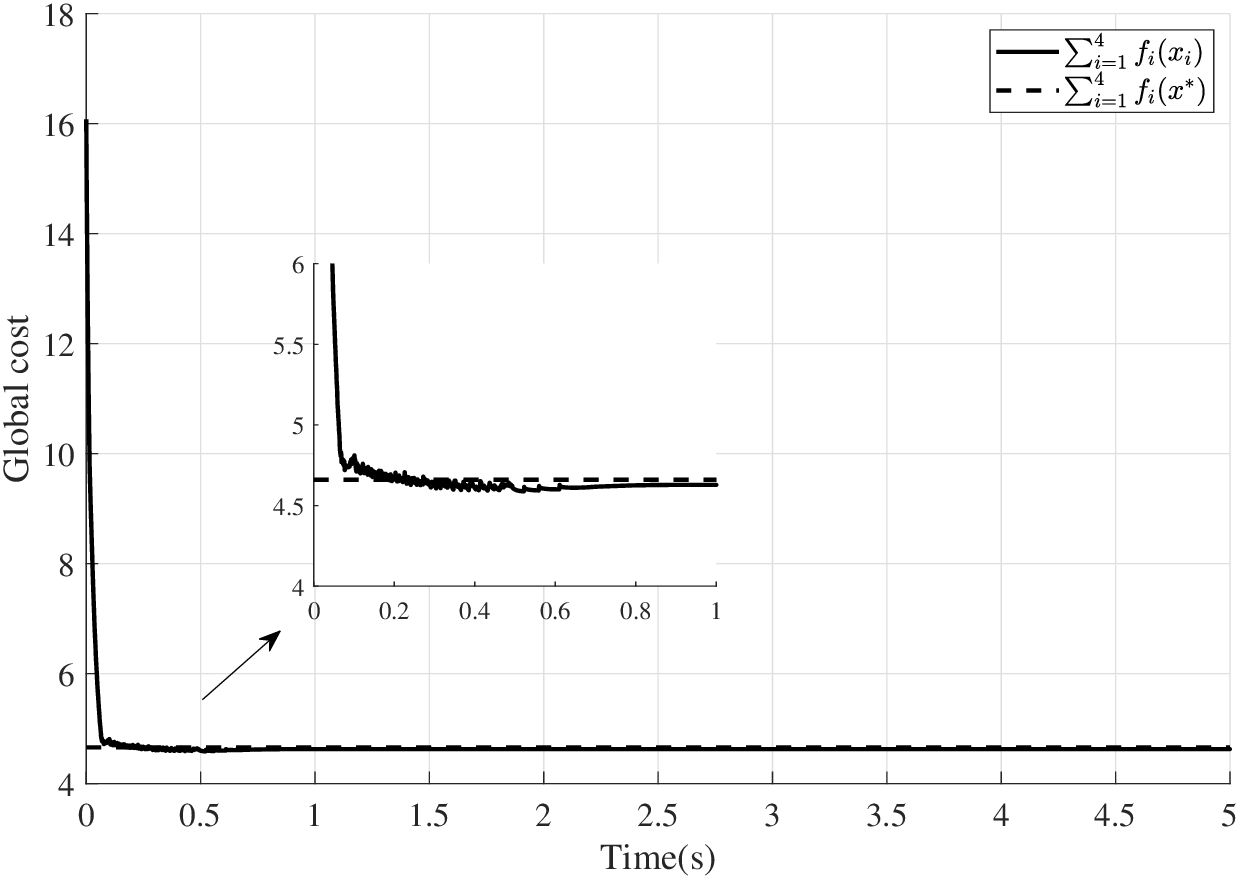}}
	\caption{Evolution of the global objective function under algorithm \eqref{Prescribed_Algorithm_time-varying}.}
	\label{B3}
\end{figure}
\begin{figure}[htb]
	\centerline{\includegraphics[width=6cm]{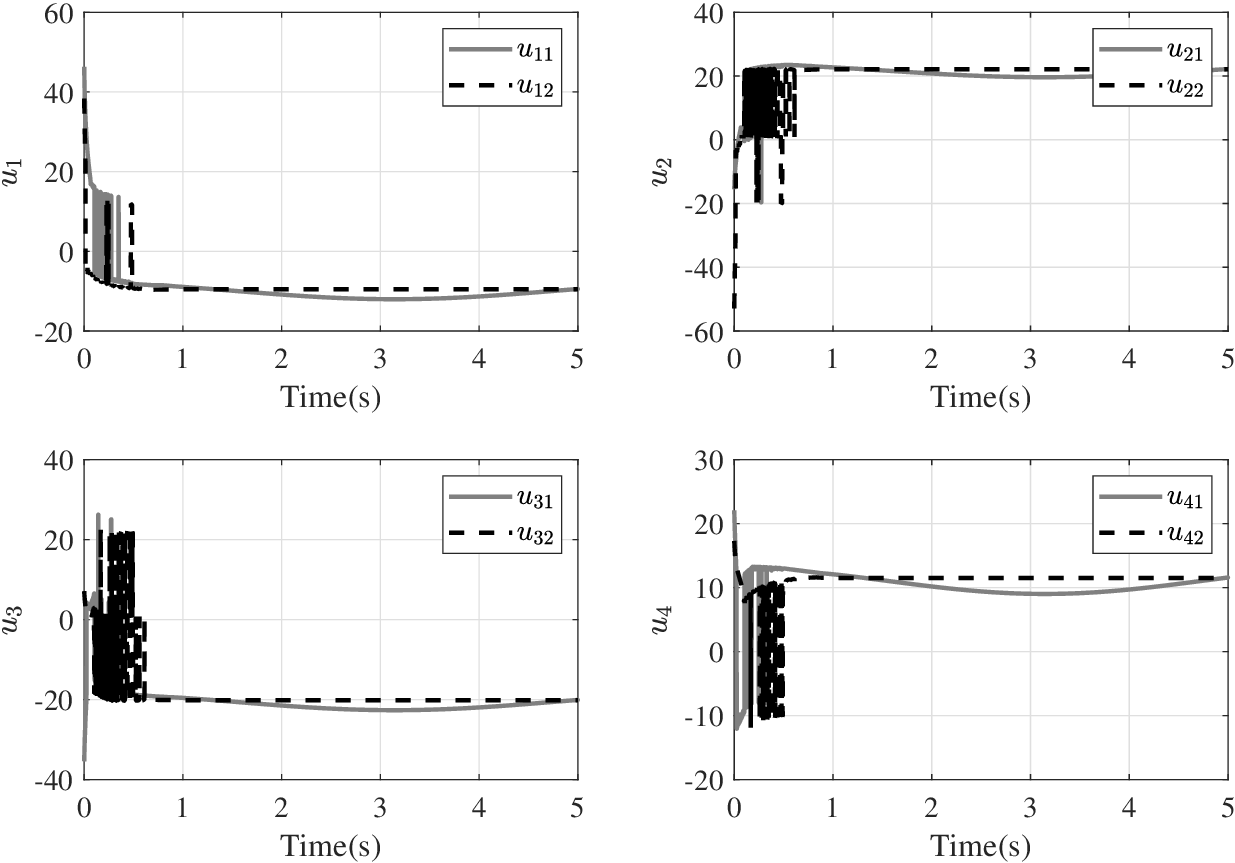}}
	\caption{Control inputs of robot $i$, $i \in \mathcal{V}$. }
	\label{B4}
\end{figure}

\subsection{Formation for target encirclement}
\begin{figure}[htb]
	\centerline{\includegraphics[width=6cm]{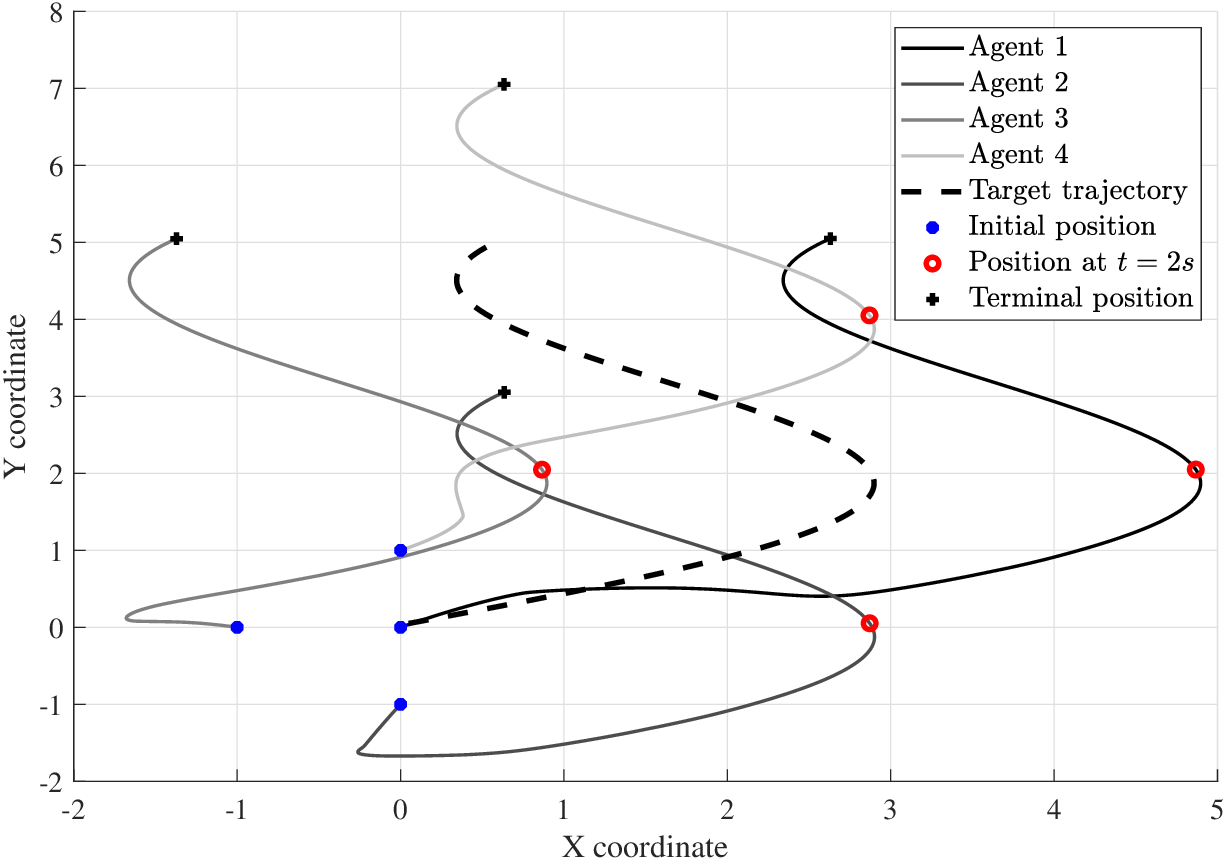}}
	\caption{Trajectories for robot $i$, $i \in \mathcal{V}$ and the target source.}
	\label{B1app}
\end{figure}

\begin{figure}[htb]
	\centering
	\subfloat[]{ \includegraphics[width=4.25cm]{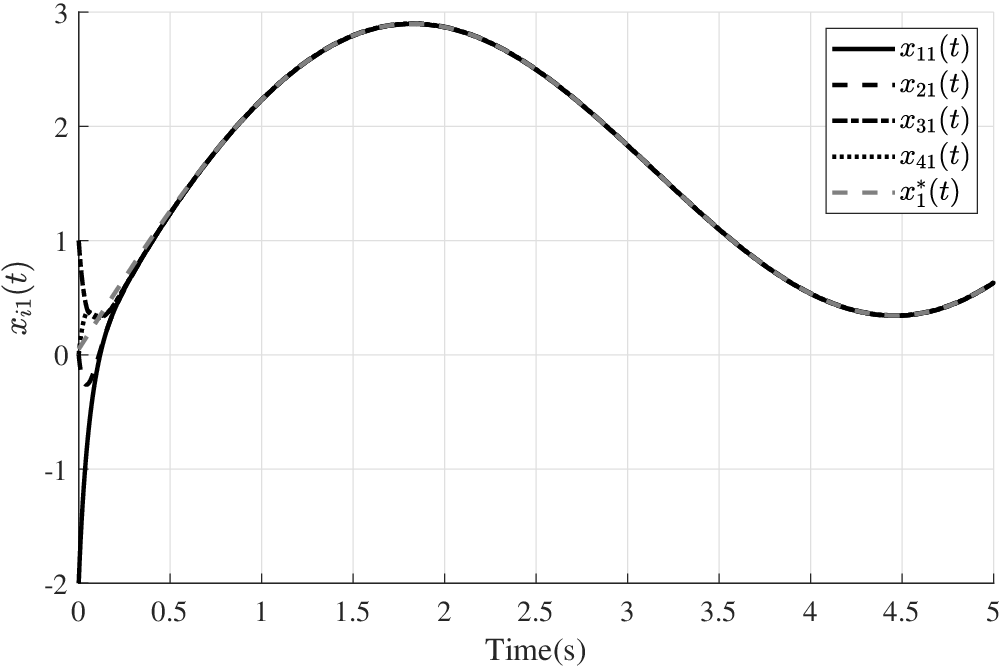} }
%	\hfil 
	\subfloat[]{ \includegraphics[width=4.2cm]{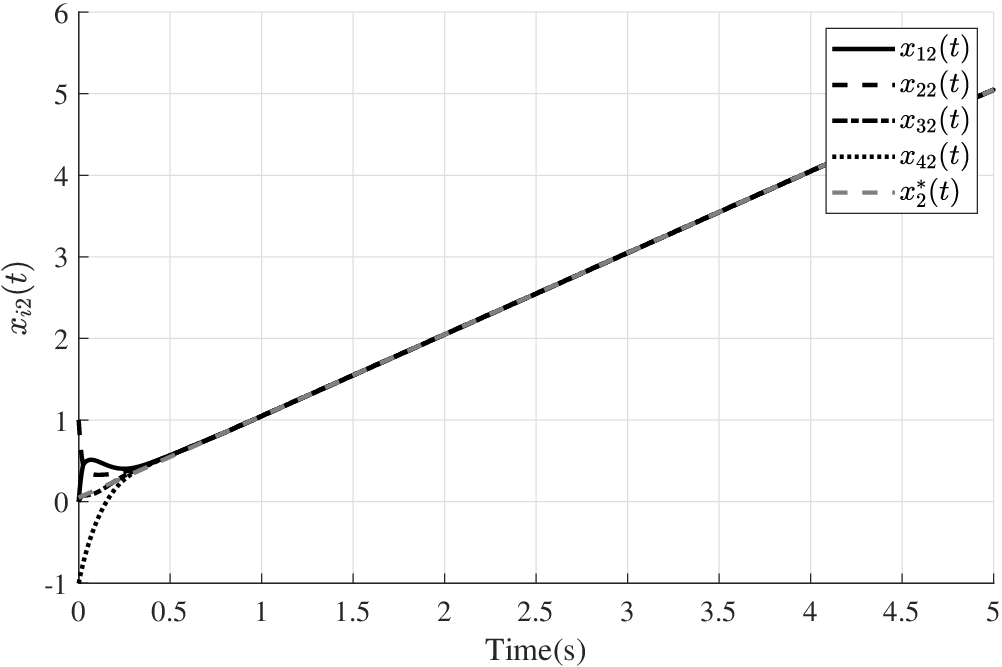} }	
	\caption{Evolution of $ x_i\left(t\right) $ under algorithm \eqref{Prescribed_Algorithm_time-varying} with boundary layer approximation. (a) $ x_{i1}\left(t\right) $. (b) $ x_{i2}\left(t\right) $ }
	\label{B2app}
\end{figure}

\begin{figure}[htb]
	\centerline{\includegraphics[width=6cm]{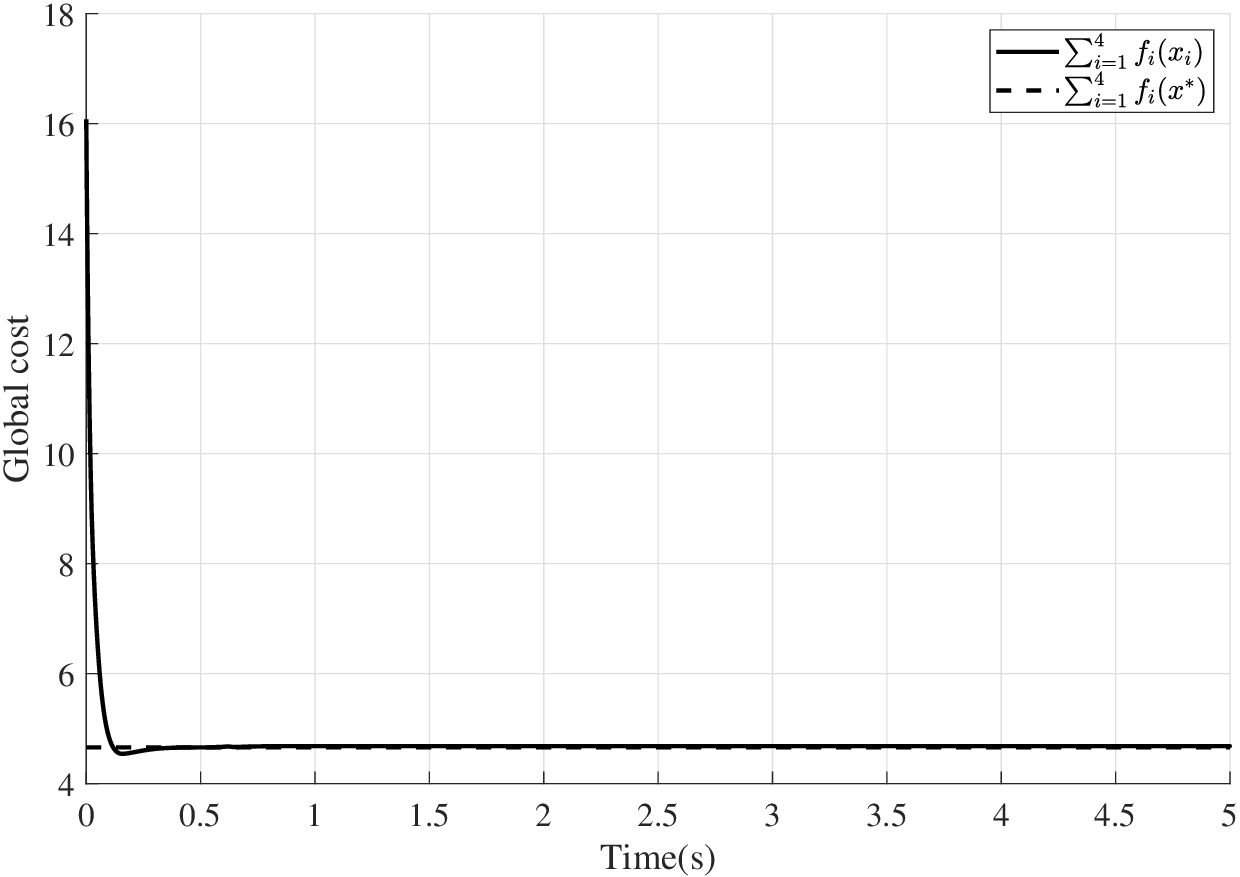}}
	\caption{Evolution of the global objective function under algorithm \eqref{Prescribed_Algorithm_time-varying} with boundary layer approximation.}
	\label{B3app}
\end{figure}

\begin{figure}[htb]
	\centerline{\includegraphics[width=6cm]{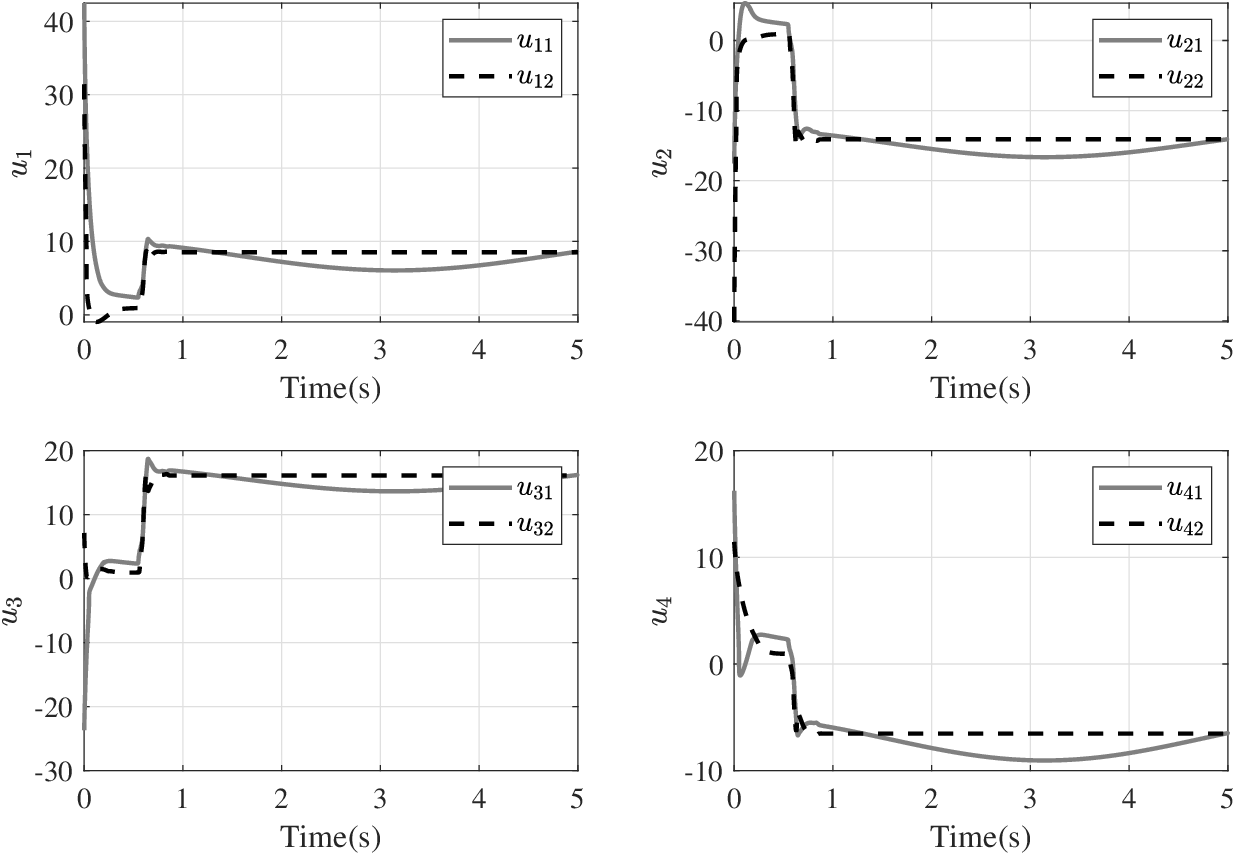}}
	\caption{Control inputs of robot $i$, $i \in \mathcal{V}$. }
	\label{B4app}
\end{figure}

A case of achieving the formation for target encirclement via distributed time-varying optimization algorithms is provided to verify the proposed method \eqref{Prescribed_Algorithm_time-varying}.
Consider a group of four robots interacting through undirected connected graph shown in Fig. \ref{fig_communication_b}. The dynamics of each robot is modeled as
$$ \dot p_i(t) = u_i(t), i \in \mathcal{V},$$ where $ p_i(t) \in \mathbb{R}^2 $ and $ u_i(t) \in \mathbb{R}^2 $ are the position and the control input of $i$th robot, and $ \mathcal{V} = \{1,2,3,4\}$. The initial positions of robots are set as $p_1(0)=[0,0]^T$, $p_2(0)=[0,-1]^T$, $p_3(0)=[-1,0]^T$ and $p_4(0)=[0,1]^T$, and the target source to be contained is defined as $p^*(t) = [2\sin (t) + 0.5t, t]^T $. The goal of the group is to achieve  $ \lim\nolimits_{t \to T_m}\left\| p_i(t)-h_i-p^*(t) \right\|^2 = 0, i \in \mathcal{V}$, where $h_i = 2[\cos (2\pi i/N),\sin (2\pi i/N)]^T $ are the formation configuration of robots. Each robot's observation of the target is $\hat{p}_i^*(t)$, which is imprecise. It is assumed that $\hat{p}_i^*(t) = p^*(t) + \varpi _i $, where $ \varpi _1 = [1,1]^T $, $\varpi _2 = [-1,-1]^T$, $\varpi _3 = [0.5,0.5]^T$  and $\varpi _4 = [-0.3,-0.3]^T$. Define $ x_i(t) = p_i(t)-h_i$ and $f_i(x_i,t) = \left\| x_i-\hat{p}_i^*(t) \right\|^2$. Note that $\dot p_i(t) = \dot x_i(t)$, so the distributed optimization algorithm \eqref{Prescribed_Algorithm_time-varying} for solving $ \min \sum_{i=1}^{4} f_i(x_i,t)$  can be treated as the control law of robot to achieve a relatively accurate enveloping configuration. The parameters are considered as $ \eta=0.5 $, $  p =0.3$, $ T_m =2 $, $ c = 3 $, and $ \mu =44$. The position trajectories of the robots and the target source are shown in Fig. \ref{B1}, from which it can be observed that all robots have achieved encirclement of the target in formation. The evolution of the optimizer and global cost for the optimization problem are depicted in Fig. \ref{B2} and \ref{B3}, respectively. Finally, the control inputs of four robots are given in Fig. \ref{B4}. Due to the discontinuity of the sign function, there is chattering phenomenon in the control input. The validity of Algorithm \eqref{Prescribed_Algorithm_time-varying} with boundary layer approximation is verified by keeping the experimental settings unchanged. The simulation results are shown in Fig. \ref{B1app}-\ref{B4app}, from which it can be observed that each robot moves to the desired position and the control input becomes significantly smoother.


\section{Conclusion}
This paper investigates the distributed optimization algorithms with predefined-time convergence and less communication requirements. Involving a sliding manifold to ensure that the sum of gradients approaches zero, a distributed algorithm is proposed to achieve global optimal consensus. The result is extended to apply to time-varying cost functions by introducing non-smooth consensus terms and local gradients prediction. The proposed algorithms only need the primal states to be communicated, and the convergence time can be specified in advance under the uniform time.

In future research, we are interested in obtaining a discrete-time implementation of the proposed algorithm with formal convergence guarantees. Predefined-time distributed optimization algorithms with state-dependent interactions are worth investigation.

\ifCLASSOPTIONcaptionsoff
  \newpage
\fi



%\bibliographystyle{IEEEtran}
\bibliographystyle{IEEEtran}
\bibliography{bare_jrnl_Final_Revised_SGN}
%\begin{thebibliography}{1}
%\end{thebibliography}


\begin{IEEEbiography}[{\includegraphics[width=1in,height=1.25in,clip,keepaspectratio]{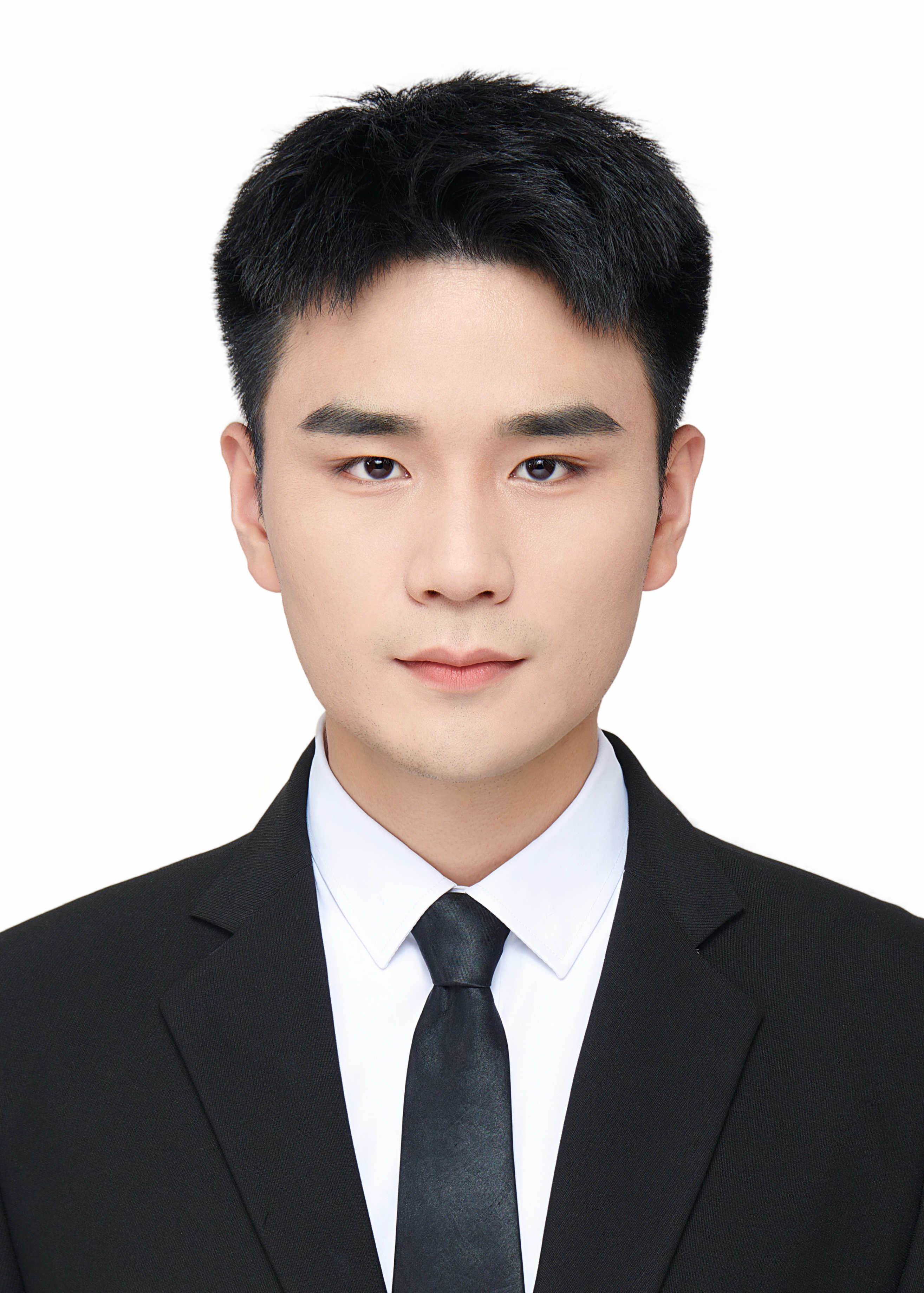}}]{Renyongkang Zhang}
	received the B.S. degree in automation from Northeastern University, Qinhuangdao, China, in 2020. He is currently pursuing the Ph.D. degree in control science and engineering at the College of Information Science and Engineering, Northeastern University, Shenyang, China. 
	
	His current research focuses on control and optimization of multi-agent systems, and with application to intelligent transportation. He was a recipient of the Best Paper Award at the \emph{5th International Conference on Industrial Artificial Intelligence (IAI 2023)}. 
	
	%	His current research focuses on control and optimization of multi-agent systems, and with application to intelligent transportation. He serves as a reviewer of several journals and conferences such as \emph{Control Engineering Practice}, \emph{Information Sciences}, \emph{ISA Transactions}, \emph{IEEE Intelligent Vehicles Symposium}, etc. He was a recipient of the Best Paper Award at the \emph{5th International Conference on Industrial Artificial Intelligence (IAI 2023)}. 
	
	%	He was a recipient of the Best Paper Award at the \emph{5th International Conference on Industrial Artificial Intelligence (IAI 2023)}. He serves as a reviewer of several journals and conferences such as \emph{Control Engineering Practice}, \emph{Information Sciences}, \emph{ISA Transactions}, \emph{IEEE Intelligent Vehicles Symposium}, etc. 
\end{IEEEbiography}

\begin{IEEEbiography}[{\includegraphics[width=1in,height=1.25in,clip,keepaspectratio]{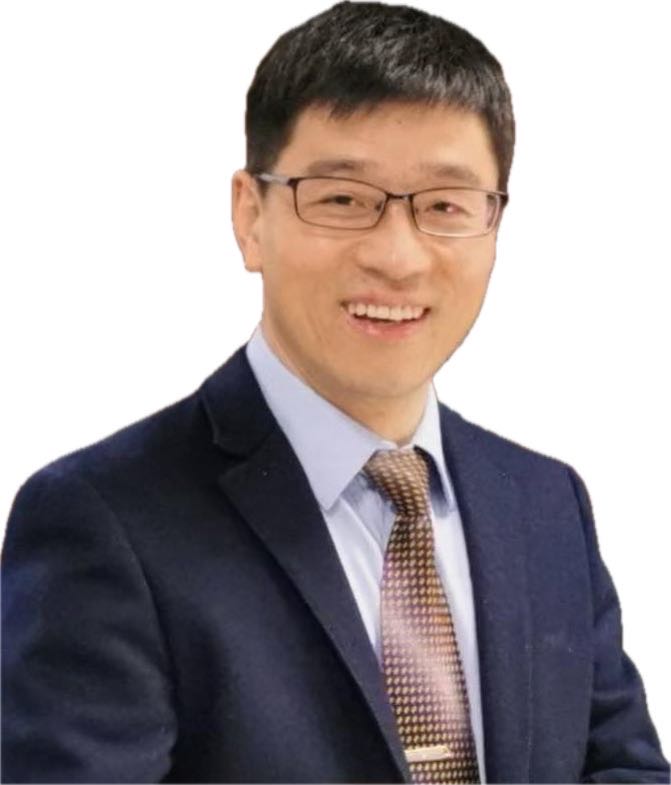}}]{Ge Guo}
	received the B.S. degree and the PhD degree from Northeastern University, Shenyang, China, in1994 and 1998, respectively.
	
	From May 2000 to April 2005, he was with Lanzhou University of Technology, China, as a Professor and the director of the Institute of Intelligent Control and Robots. He then joined Dalian Maritime University, China, as a Professor with the Department of Automation. Since 2018, he has been a Professor with Northeastern University (NEU) and  the director of the Center of Intelligent Transportation Systems of NEU. He has published over 200 international journal papers within his areas of interest, which include intelligent transportation systems, cyber-physical systems, etc. 
	
	Dr. Guo is an Associate Editor of the \emph{IEEE Transactions on Intelligent Transportation Systems}, the \emph{IEEE Transactions on Vehicular Technology}, the \emph{IEEE Transactions on Intelligent Vehicles}, the \emph{Information Sciences}, the \emph{IEEE Intelligent Transportation Systems Magazine}, the \emph{ACTA Automatica Sinica}, the \emph{China Journal of Highway and Transport} and the \emph{Journal of Control and Decision}. He won a series of awards including the CAA Young Scientist Award, the First Prize of Natural Science Award of Hebei Province, the first Prize of Science and Technology Progress Award of Gansu Province, etc.
\end{IEEEbiography}

\begin{IEEEbiography}[{\includegraphics[width=1in,height=1.25in,clip,keepaspectratio]{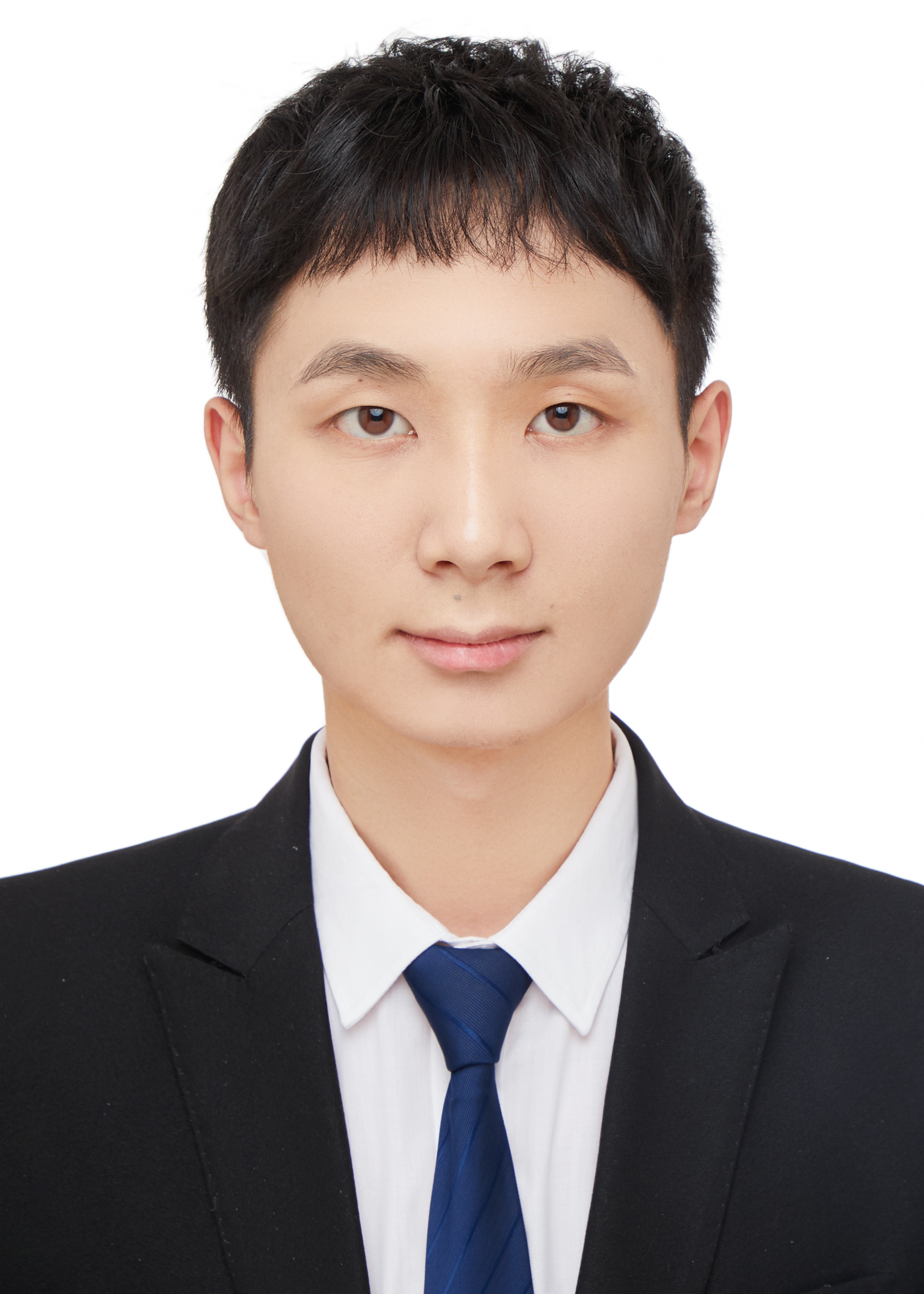}}]{Zeng-Di Zhou}
	received the B.S. degree in measurement and control technology and instrument from Northeastern University, Qinhuangdao, China, in 2021. He is currently pursuing the Ph.D. degree in control science and engineering at the College of Information Science and Engineering, Northeastern University, Shenyang, China. 
	
	His current research focuses on control and optimization of multi-agent systems. 
\end{IEEEbiography}

%\newpage


%\vfill


\end{document}